\documentclass[ejs]{imsart}

\RequirePackage[numbers,sort&compress]{natbib}
\RequirePackage[colorlinks,citecolor=blue,urlcolor=blue]{hyperref}

\doi{10.1214/154957804100000000}
\pubyear{0000}
\volume{0}
\firstpage{1}
\lastpage{1}

\startlocaldefs

\usepackage{amsmath, amsthm, amssymb, amscd, amsfonts, bm, bbm, mathrsfs}   
\usepackage{mathtools}
\usepackage{enumerate}
\usepackage[dvipsnames]{xcolor}
\usepackage[shortlabels]{enumitem}
\usepackage[capitalize,nameinlink,noabbrev]{cleveref} 
\usepackage{ulem}

\theoremstyle{plain}
\newtheorem{thm}{Theorem}
\newtheorem{cor}{Corollary}

\newtheorem{prop}{Proposition}
\newtheorem{rem}{Remark}
\newtheorem{eg}{Example}

\newcommand{\tr}[1]{\text{tr}\left(#1\right)} 

\makeatletter
\newcommand*{\rom}[1]{\expandafter\@slowromancap\romannumeral #1@}
\newcommand{\HUGE}{\@setfontsize\Huge{40}{50}}   
\makeatother

\makeatletter
\newcommand{\labeltext}[3][]{%
	\@bsphack%
	\csname phantomsection\endcsname
	\def\tst{#1}%
	\def\refmarkup{}%
	\ifx\tst\empty\def\@currentlabel{\refmarkup{#2}}{\label{#3}}%
	\else\def\@currentlabel{\refmarkup{#1}}{\label{#3}}\fi%
	\@esphack%
	\labelmarkup{#2}
}
\makeatother

\makeatletter
\newcommand{\bianca}{\renewcommand\NAT@open{[}\renewcommand\NAT@close{]}}
\makeatother

\newcommand{\iid}{\overset{\mathsf{iid}}{\sim}} 
\newcommand{\ed}{\overset{\mathsf{d}}{=}} 

\newcommand{\pr}{\mathsf{P}}
\newcommand{\eo}{\mathsf{E}}

\newcommand{\var}{\mathsf{var}}

\newcommand{\nd}{\mathsf{N}}

\newcommand{\ap}{\alpha} 
\newcommand{\g}{\gamma} 
\newcommand{\ga}{\Gamma} 
\newcommand{\dt}{\delta} 
\newcommand{\Dt}{\Delta} 
\newcommand{\e}{\varepsilon} 
\newcommand{\ka}{\kappa} 
\newcommand{\s}{\sigma} 
\newcommand{\ld}{\lambda} 
\newcommand{\Ld}{\Lambda} 

\newcommand{\E}{\mathbb{E}} 
\newcommand{\HH}{\mathbb{H}} 
\newcommand{\N}{\mathbb{N}} 
\newcommand{\R}{\mathbb{R}} 

\newcommand{\dd}{\mathcal{D}}	
\newcommand{\ii}{\mathcal{I}}	
\newcommand{\kll}{\mathcal{L}}	
\newcommand{\xx}{\mathcal{X}}	

\endlocaldefs

\begin{document}

\newcommand{\thmautorefname}{Theorem}
\newcommand{\defnautorefname}{Definition}
\newcommand{\propautorefname}{Proposition}
\newcommand{\corautorefname}{Corollary}
\newcommand{\lemautorefname}{Lemma}
\newcommand{\remautorefname}{Remark}
\newcommand{\egautorefname}{Example}

\renewcommand{\sectionautorefname}{Section}
\renewcommand{\subsectionautorefname}{Section}
\renewcommand{\subsubsectionautorefname}{Section}
\renewcommand{\eqref}[1]{(\ref{#1})}

\allowdisplaybreaks

\begin{frontmatter}


\title{
	Bootstrap inference in functional linear regression models with scalar response under heteroscedasticity
}
\runtitle{Bootstrap in FLRM under heteroscedasticity}


\author[A]{\fnms{Hyemin}~\snm{Yeon}\corref{}\ead[label=e1]{hyeon@ncsu.edu}},
\author[B]{\fnms{Xiongtao}~\snm{Dai}\ead[label=e2]{xdai@berkeley.edu}}
\and
\author[B]{\fnms{Daniel J.}~\snm{Nordman}\ead[label=e3]{dnordman@iastate.edu}}

\address[A]{Department of Statistics, North Carolina State University\printead[presep={,\ }]{e1}}

\address[B]{Division of Biostatistics, University of California, Berkeley\printead[presep={,\ }]{e2}}

\address[C]{Department of Statistics, Iowa State University\printead[presep={,\ }]{e3}}

\runauthor{H. Yeon, X. Dai, and D. J. Nordman}

\begin{abstract}
	Inference for functional linear models in the presence of heteroscedastic errors has received insufficient attention given its practical importance;  in fact, even a central limit theorem has not been studied in this case. 
	At issue, conditional mean estimates have complicated sampling distributions due to the infinite dimensional regressors, where truncation bias and scaling issues are compounded by non-constant variance under heteroscedasticity. 
	As a foundation for distributional inference, we establish a central limit theorem for the estimated conditional mean under general dependent errors, 
	and subsequently we develop a paired bootstrap method to provide better approximations of  sampling distributions.
	The proposed paired bootstrap does not follow the standard bootstrap algorithm for finite dimensional regressors, as
	this version fails outside of a narrow window for implementation with functional regressors. 
	The reason owes to a bias with functional regressors in a naive bootstrap construction. 
	Our bootstrap proposal incorporates debiasing and thereby attains much broader validity and flexibility with truncation parameters for inference under heteroscedasticity; 
	even when the naive approach may be valid, the proposed bootstrap method performs better numerically. 
	The bootstrap is applied to construct confidence intervals for centered projections and for conducting hypothesis tests for the multiple conditional means. 
	Our theoretical results on bootstrap consistency are demonstrated through simulation studies and also illustrated with a real data example.
\end{abstract}

\begin{keyword}[class=MSC]
	\kwd[Primary ]{62R10}
	\kwd{62G09}
	\kwd[; secondary ]{62E20}
\end{keyword}

\begin{keyword}
	\kwd{Asymptotic normality}
	\kwd{Bias correction}
	\kwd{Bootstrapping pairs}
	\kwd{Functional data analysis}\kwd{Multiple testing}
	\kwd{Scalar-on-function regression}
\end{keyword}


\tableofcontents

\end{frontmatter}


\newcommand{\tred}[1]{{\color{red}#1}} 
\newcommand{\tblue}[1]{{\color{blue}#1}} 
\newcommand{\tgreen}[1]{{\color{ForestGreen}#1}} 


\section{Introduction} \label{sec1}

In classical linear models, bootstrap methods have been developed for several decades 
under either homoscedastic or heteroscedastic error assumptions.
  Residual and paired bootstrap methods 
were originally studied by \cite{free81} for approximating the   distribution of the least square estimator in multiple linear regression models. 
These bootstraps are intended, respectively, for handling homoscedastic or heteroscedastic error cases. 
Both bootstraps
have been investigated in other contexts as well,
such as nonparametric \citep{HM91} or high-dimensional \citep{DBZ17} regression problems.
In a functional linear regression model (FLRM), 
bootstrap inference is likewise valuable but also   complicated 
due to the infinite dimensionality of the underlying function space. 
A main issue with functional regressors is that a truncation bias arises in estimators of the conditional mean, due to the infinite dimensional regressors and slope functions involved, which imposes challenges for even central limit theorems \citep{CMS07, YDN23RB}.  
Existing works on both the central limit theorem (CLT) and (residual/wild) bootstrap for functional linear regression models (FLRMs) have focused exclusively on  homogeneous error variance models \citep{GM11, YDN23RB}, 
while either avoiding or accommodating   bias issues. 
In fact, beyond homoscedasticity, more stringent conditions of independence between regressors and errors are also commonly imposed in FLRM literature \citep{CH06, cai:12, HH07} and especially for hypothesis testing  \citep{CFMS03, CGS04, HMV13, lei14,LL22}.
However,  heteroscedastic error variances are commonly observed in practice.

\begin{figure}[b!]
	\centering
	\includegraphics[width=0.6\linewidth]{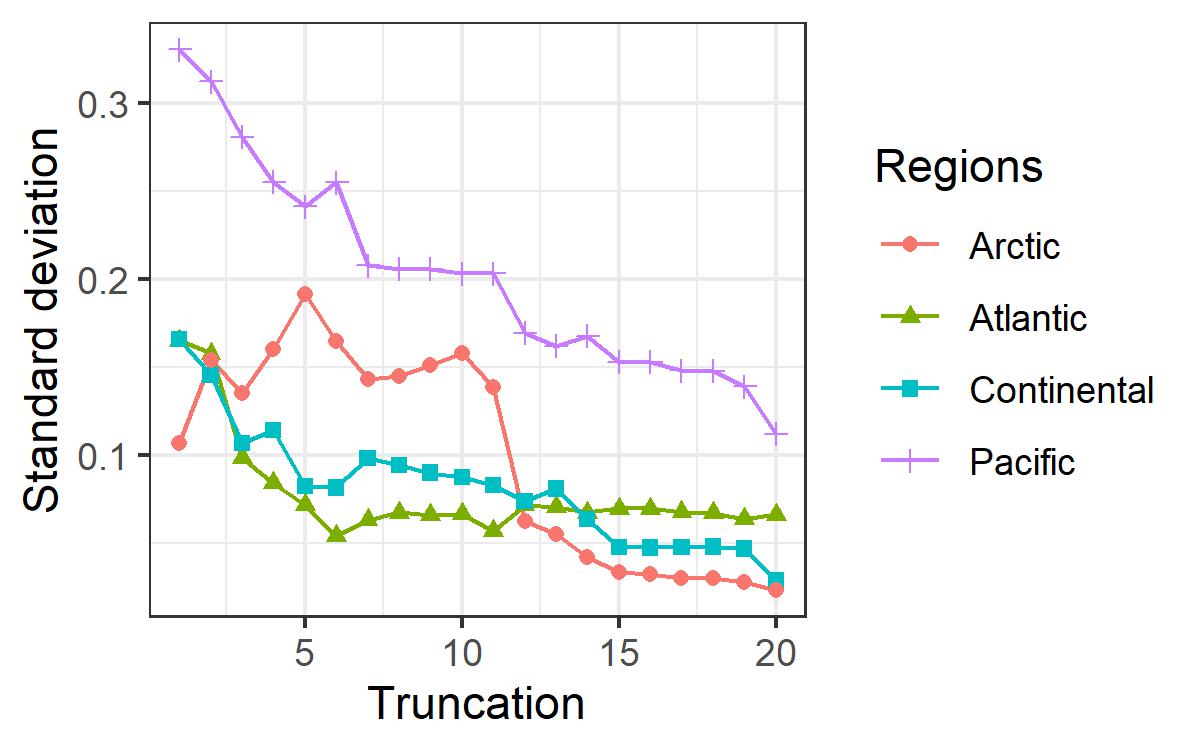}
	\caption{Estimated standard deviations for each region over different truncation levels used in estimation}
	\label{fig_rda_cw_estSD}
\end{figure}

For illustration, \autoref{fig_rda_cw_estSD} shows the estimated standard deviations 
of residuals from a  FLRM fit to a Canadian weather dataset (cf.~\autoref{sec5}) over different geographical regions. 
Each regressor curve represents averaged daily temperatures measured at a different location contained in one of the four regions in Canada: Atlantic, Continental, Pacific, and Arctic regions,
where the associated response is the total annual precipitation on the log scale. 
As variances appear to differ across regions, it seems natural here to avoid homoscedastic error models.

To the best of our knowledge, heteroscedastic error conditions have not received much formal consideration in the FLRM literature, with perhaps the exception of work on weighted least squares by \cite{DHT09}, which does not discuss distributional inference.  For example, while a CLT for projection estimates is again available for FLRMs in the homoscedastic case \citep{CMS07, KH16a, YDN23RB},
a counterpart foundational result does not yet exist under heteroscedasticity.    
One might further  anticipate that previous bootstrap theory  under homoscedasticity  does not directly apply for the inference in FLRMs under heteroscedasticity.  We show this to be the case, 
which necessitates our new development of a CLT and resampling theorems.  As in the homoscedastic setting, resampling approximations in FLRMs are   valuable under heteroscedasticity for capturing complicated sampling distributions of mean estimators, as current bootstraps from  the  homoscedastic case  become invalid \citep{GM11,   KH16b, YDN23RB}.

To bootstrap FLRMs in the presence of heteroscedastic errors, a paired bootstrap method can be considered, 
similar to the paired bootstrap for usual multiple linear regression models \citep{free81}.
Paired bootstrap has indeed been applied for  different functional regression models  
\citep{MCL08, RO10, WH17},
though without any theoretical development or justification.  
This latter point is important, because we show here that, surprisingly, a naive/standard implementation of paired bootstrap, adopted directly from the usual multiple regression case, can fail to provide valid inference for mean estimates under FLRMs if the truncation parameters are not set appropriately in a certain narrow and restricted way,
in contrast to the case
of finite-dimensional multiple linear regression \citep{free81}. 
In fact, as sample sizes increase 
the   distance   between   naive paired bootstrap   and  true sampling distributions  may not converge to zero  as typically expected, but rather to a random number, unless associated tuning parameters are set in a specific manner.
The problem arises from a construction bias in the bootstrap world with FLRMs which relates to, but is a separate issue from, the truncation bias inherent to 
the functional principal component regression estimator $\hat{\beta}_{h_n}$ 
of the slope (cf.~\autoref{sec2}).   
Such failure of bootstrap due to bias issues has been observed in other bootstrap works with complicated regressions, such as nonparametric \citep{HM91, HM93, zie01}, quantile \citep{WMQ18}, penalized linear \citep{CL10, CL11, camp15}, and high-dimensional linear \citep{DBZ17} regression models, 
though the approaches of correcting bootstrap bias can differ. 
In some problems, the extent of the bias in paired bootstrap is such that this bootstrap must be discarded (cf.~\cite{HM91, HM93}).
This motivates our new development of paired bootstrap for FLRMs with heteroscedastic errors, which remedies the bias problem by modifying a bootstrap estimating equation to define a bootstrap estimator.

Under a general heteroscedastic error assumption, we study asymptotic  inference in FLRMs with scalar response. 
In particular, we first establish a CLT for the estimator of the conditional mean $\mu(X_0) \equiv \ap + \langle \beta, X_0 \rangle$
 with $X_0$ being a (random) new regressor function.
This provides a foundation for our bootstrap results and, more broadly, supports inference for FLRMs under dependent errors.
Our main bootstrap result is to develop a modified paired bootstrap to approximate the sampling distribution of mean estimators based on the functional principal component regression (FPCR) estimator $\hat{\beta}_{h_n}$  of the slope function $\beta$ \citep{CFS99, CH06, HH07, CMS07, GM11, KH16a, KH16b, YDN23RB}, where $h_n$ denotes a truncation level involved in the estimation procedure.
For flexibility and   better practical performance, the modified paired bootstrap incorporates an important debiasing step.
In heteroscedastic cases, our numerical studies suggest that the paired bootstrap performs better than the residual bootstrap and normal approximation, while also maintaining good coverage  in homoscedastic cases.
The proposed paired bootstrap also numerically outperforms the naive version even when the latter is appropriately tuned.

As an application of the paired bootstrap method, 
we also treat a testing problem about evaluating the constancy of conditional means over a collection of target functional regressors; 
this assesses a null hypothesis that the slope function   has no effect on (or, in a sense, is orthogonal to) some specified regressors. 
In this problem, the bootstrap combines several simultaneous estimation steps into one test, which would otherwise be distributionally intractable through normal approximations.
The bootstrap  also has the advantage of enforcing the null hypothesis in re-creating a reference distribution for testing, which can be useful for controlling size and boosting power.   
Our development in this testing problem is distinguishable from the previous works on hypothesis testing in FLRMs  \citep{CFMS03, CGS04, HMV13, lei14,LL22}: 
the latter tests are limited to independent error scenarios and often restrict claims to global nullity $\beta=0$  based on  certain basis functions.  
In contrast, our testing method also allows  hypotheses about $\beta$ to be defined   with more arbitrarily specified functional regressors along with dependence between regressors and errors.
This is useful in practice for assessing 
how  projections of $\beta$ (or conditional means) may differ   
as predictor levels are varied, which may not be  addressable by a global test of $\beta$.

\autoref{sec2} provides background  on  FLRMs and the paired bootstrap method under heteroscedasticity. 
\autoref{ssec_2_3} presents the main distributional results regarding  regarding the consistency of the paired bootstrap method and the failure of the naive bootstrap.  With suitable scaling, a general CLT is also established for  conditional  mean estimators.
We then give a bootstrap procedure in \autoref{sec3} for hypothesis testing  of multiple conditional means.
Numerical results are provided in \autoref{sec4},  
while \autoref{sec5} illustrates the paired bootstrap method with a real dataset having potential  heteroscedasticity
Concluding remarks are offered in \autoref{sec6}.
Some proofs for the main results are given in Appendix, while further details of the proofs and extended numerical results can be found in a supplement \cite{supp}.  
An R package \cite{pack} is provided to construct confidence intervals for FLRM response means and test multiple projections   of $\beta$ based on paired bootstrap.

\section{Description of FLRMs and  bootstrap} \label{sec2}

We start with the description of functional linear regression models (FLRMs) under heteroscedastic error variances in \autoref{ssec_2_1}, and the paired bootstrap for estimated conditional means appears in \autoref{ssec_2_2}.

\subsection{FLRMs under heteroscedasticity} \label{ssec_2_1}

Consider the following FLRM
\begin{align} \label{eq_model_intercept}
	Y = {\ap} + \langle \beta, X \rangle + \e,
\end{align}
where $Y$ is a scalar-valued response; $X$ is a regressor function taking values in a separable Hilbert space $\HH$ with inner product $\langle \cdot, \cdot \rangle$; {$\alpha$ is the intercept;} and $\beta \in \HH$ is the slope function.
The error term $\e$ has $\eo[\e|X] = 0$ but its distribution can otherwise depend on $X$; for example, heterogeneous conditional variances of the error $\e$  given the regressor $X$ is allowed, that is, $\s^2(X) \equiv \eo[\e^2|X]$ may depend on the regressor $X$. 

Define the tensor product $x \otimes y:\HH \times \HH \to \HH$ between two elements $x,y \in \HH$ as a bounded linear operator $z \mapsto (x \otimes y)(z) = \langle z, x \rangle y$ for $z \in \HH$ and denote $x^{\otimes 2} \equiv x \otimes x$ for $x \in \HH$. 
Under the assumption $\eo[\|X\|^2]<\infty$, where $\|\cdot\|$ is the induced norm in $\HH$, the covariance operator $\ga \equiv \eo[(X - \eo[X])^{\otimes 2}]$ is self-adjoint, non-negative definite, and Hilbert--Schmidt, and hence, compact (cf.~\cite{HE15}). Then, $\ga$ admits the following spectral decomposition 
\begin{align*}
	\ga = \sum_{j=1}^\infty \g_j \pi_j
\end{align*}
with $\pi_j \equiv \phi_j \otimes \phi_j$,
where $\g_j$  and $\phi_j$ are the $j$-th eigenvalue and eigenfunction of $\ga$ for $j = 1,2,\dots$. 
Here, the set $\{\phi_j\}$ of eigenfunctions is an orthonormal system of $\HH$ and $\{\g_j\}$ is a non-negative non-increasing sequence with $\g_j \to 0$ as $j\to\infty$.
The functional version of normal equations can be written as 
\begin{align}
	\Dt = \ga \beta
\end{align}
from the model \eqref{eq_model_intercept}, where $\Dt \equiv \eo[(Y - \eo[Y])(X-\eo[X])]$ is the cross-covariance function between $X$ and $Y$.
Under the model identifiability assumption $\ker \ga = \{0\}$ \citep{CFS99, CFS03,CMS07} (see Condition~\ref{condA1_model} of \autoref{sssec_2_3_1}), the slope function is then given as
\begin{align*}
	\beta = \ga^{-1} \Dt.
\end{align*}

The functional principal component regression 
(FPCR) estimator of  $\beta$ has been widely studied in the literature \citep{ CH06,  CFS99, CMS07, HH07}. To define the estimator, we suppose that the data pairs $\{(X_i, Y_i)\}_{i=1}^n$ are independently and identically distributed under the FLRM \eqref{eq_model_intercept}, that is,
\begin{align} \label{eq_model_obs}
	Y_i = \ap + \langle \beta, X_i \rangle + \e_i, \quad i=1, \dots, n.
\end{align}
The sample versions of $\ga$ and $\Dt$ are defined as $\hat{\ga}_n \equiv n^{-1} \sum_{i=1}^n (X_i - \bar{X})^{\otimes 2}$ and $\hat{\Dt}_n \equiv n^{-1} \sum_{i=1}^n (Y_i - \bar{Y}) (X_i - \bar{X})$, where $\bar{X} \equiv n^{-1} \sum_{i=1}^n X_i$, $\bar{Y} \equiv n^{-1} \sum_{i=1}^n Y_i$, and $x^{\otimes 2} \equiv x \otimes x$ for $x \in \HH$.
The sample covariance operator $\hat{\ga}_n$ also admits spectral decomposition $\hat{\ga}_n = \sum_{j=1}^n \hat{\g}_j \hat{\pi}_j$ with $\hat{\pi}_j \equiv \hat{\phi}_j \otimes \hat{\phi}_j$, where $\hat{\g}_j \geq 0$ is the $j$-th sample eigenvalue and $\hat{\phi}_j \in \HH$ is the corresponding eigenfunction. By regularizing the inversion of $\hat{\ga}_n$, the FPCR estimator of $\beta$ is defined as 
\begin{align} \label{eq_fpcr}
	\hat{\beta}_{h_n} \equiv \hat{\ga}_{h_n}^{-1} \hat{\Dt}_n
\end{align}
where $\hat{\ga}_{h_n}^{-1} \equiv \sum_{j=1}^{h_n} \hat{\g}_j^{-1} \hat{\pi}_j$ is a finite approximation of $\ga^{-1} \equiv \sum_{j=1}^\infty \g_j^{-1} \pi_j$,
{
	and the intercept is estimated by $\hat{\ap}_{h_n} \equiv \bar{Y} - \langle \hat{\beta}_{h_n}, \bar{X} \rangle$.
	}
Here, $h_n$ is the number of eigenpairs used in estimation, which represents a truncation level 
\cite{ CH06,  CFS99, CMS07, HH07}.

\subsection{Paired bootstrap procedure} \label{ssec_2_2}

For FLRMs with homoscedastic errors, the residual bootstrap is natural \citep{GM11, YDN23RB}, where this bootstrap re-creates data, e.g., $Y_i^* = \hat{\ap}_{h_n} + \langle  \hat{\beta}_{h_n} , X_i \rangle + \varepsilon_i^*$, 
through bootstrap error terms $\e^*$ as independent draws from an appropriate set of residuals.      
However, under heteroscedastic errors, a different bootstrap approach is necessary, akin to the regression case with Euclidean vectors \citep{free81}.
Similar to that setting for capturing response variances that may differ conditionally over regressors, we consider a paired bootstrap (PB) method for inference in FLRMs. 
To the best of our knowledge, the theory for PB in FLRMs has been studied only once by \cite{GGMG12}, but their application does not consider  slope functions
and the errors therein are homoscedastic.
For estimating means or projections under the FLRM with heteroscedastic errors, we explain next how the PB generally requires careful consideration in order to be valid.

To implement the PB, we define bootstrap data pairs $\{(X_i^*, Y_i^*)\}_{i=1}^n$ by uniform draws from the original data $\{(X_i, Y_i)\}_{i=1}^n$ with replacement. The bootstrap counterparts of sample moments are then given as $\hat{\ga}_n^* \equiv n^{-1} \sum_{i=1}^n (X_i^* - \bar{X}^*)^{\otimes 2}$ and $\hat{\Dt}_n^* \equiv n^{-1} \sum_{i=1}^n (Y_i^* - \bar{Y}^*) (X_i^* - \bar{X}^*)$ where $\bar{X}^* \equiv n^{-1} \sum_{i=1}^n X_i^*$ and $\bar{Y}^* \equiv n^{-1} \sum_{i=1}^n Y_i^*$.
From the spectral decomposition of $\hat{\ga}_n^*$, we define a regularized inverse of $\hat{\ga}_n^*$ with trunction level $h_n$ as
\begin{align*}
	(\hat{\ga}_{h_n}^*)^{-1} \equiv \sum_{j=1}^{h_n} (\hat{\g}_j^*)^{-1} (\hat{\phi}_j^* \otimes \hat{\phi}_j^*),
\end{align*}
where $\hat{\g}_j^*$ and $\hat{\phi}_j^*$ are the $j$-th eigenvalue and the corresponding eigenfunction of $\hat{\ga}_n^*$.  This represents a direct bootstrap analog of $\hat{\ga}_{h_n}^{-1}$ in \eqref{eq_fpcr}.

An initial, though naive, bootstrap version  $\hat{\beta}_{h_n,  naive}^*$ of the FPCR estimator $\hat{\beta}_{h_n}$ can be found  as
\begin{align}
	\hat{\beta}_{h_n,  naive}^* \equiv (\hat{\ga}_{h_n}^*)^{-1} \hat{\Dt}_n^* \label{eq_fpcr_naive}
\end{align}
by
directly imitating the  definition of $\hat{\beta}_{h_n}$ in \eqref{eq_fpcr} with bootstrap data.  The validity of this naive bootstrap, though, requires caution.  The issue is that, in the bootstrap world, we require a 
bootstrap version $\beta^*$ of the true parameter $\beta$
and, for flexibility, one might consider a FPCR estimator $\beta^*\equiv \hat{\beta}_{g_n} \equiv \hat{\ga}_{g_n}^{-1}  \hat{\Dt}_n$ determined by a general  
truncation level $g_n$ in  \eqref{eq_fpcr}.  It turns out that  the naive  bootstrap estimator $\hat{\beta}_{h_n, naive}^*$ must   restrictively use a bootstrap parameter   $\beta^*\equiv \hat{\beta}_{g_n}$ defined by $g_n=h_n$ due to    a construction bias in the naive bootstrap.  Unless the bootstrap parameter $\beta^*$ 
is specifically chosen as $\hat{\beta}_{h_n}$, which imposes limitations for implementation and numerical performance, the naive bootstrap   will be provably biased with
an adverse effect on inference 	
(cf.~\autoref{prop_pb_naive_fail}).

In order to define a more versatile bootstrap version $\hat{\beta}_{h_n}^*$ of the FPCR estimator $\hat{\beta}_{h_n}$, 
let $\hat{\beta}_{g_n}$ again denote a FPCR estimator, similar to   $\hat{\beta}_{h_n}$ from \eqref{eq_fpcr} but based on a truncation $g_n$ rather than $h_n$, in order 
   to play the role $\beta^*$ of the slope function $\beta$ in the bootstrap world.  The level of truncation $g_n$ used in a bootstrap version $\beta^*=\hat{\beta}_{g_n}$ of $\beta$ becomes a consideration  because  $\beta$ is infinite dimensional while any FPCR estimator $\hat{\beta}_{g_n}$ is finite-dimensional.    
While possible to choose $g_n=h_n$, more flexibility with $\beta^*=\hat{\beta}_{g_n}$ for $g_n$ smaller than $h_n$ can later provide   better PB approximations.   
For any estimator $\hat{\beta}_{g_n}$ version used to mimic $\beta$, 
the PB analog $\hat{\beta}_{h_n}^*$ of the original-data estimator $\hat{\beta}_{h_n}$ needs to be appropriately defined to avoid a construction bias in resampling. 
To correct this  bias, we adapt a modification of \cite{shor82} for defining bootstrap M-estimators  through  adjusted bootstrap estimating equations;
see also \cite{KL94, lahiri92} or \cite[Section~4.3]{lahiri03}.


To  define a modified bootstrap version $\hat{\beta}_{h_n}^*$ of the FPCR estimator  $\hat{\beta}_{h_n}$, we first observe that the slope function $\beta = \ga^{-1} \Dt$ can be prescribed as the solution 
to the estimating equation $\eo[\Psi_i(\beta; \mu_X, \mu_Y)] = 0$, where 
\begin{equation}
	\label{eqn:mest1}
	\Psi_i(\beta; \mu_X, \mu_Y) \equiv (X_i - \mu_X) (Y_i - \mu_Y) - (X_i - \mu_X)^{\otimes 2} \beta
\end{equation}
is an estimating function with $\mu_X \equiv \eo[X]$ and $\mu_Y \equiv \eo[Y]$.
A direct bootstrap counterpart of this estimating function is given by, say, 
\begin{align*}
	\check{\Psi}_i^*(\beta; \bar{X}, \bar{Y}) 
	\equiv (X_i^* - \bar{X}) (Y_i^* - \bar{Y}) - (X_i^* - \bar{X})^{\otimes 2} \beta,
\end{align*}
where $\bar{X} \equiv n^{-1} \sum_{i=1}^n X_i$ and $\bar{Y} \equiv n^{-1} \sum_{i=1}^n Y_i$.
A key observation is that, while $\beta = \ga^{-1} \Dt$ is the solution to the equation $\eo[\Psi_i(\beta; \mu_X, \mu_Y)] = 0$, an estimator $\hat{\beta}_{g_n} \equiv \hat{\ga}_{g_n}^{-1} \hat{\Dt}_n$, playing the role of $\beta$ in the bootstrap world, will not generally be a solution to the equation 
\begin{align*}
	\hat{\Dt}_n - \hat{\ga}_n \beta \equiv \eo^*[\check{\Psi}_i^*(\beta; \bar{X}, \bar{Y})] =  0
\end{align*} 
due to the finite dimensionality of $\hat{\beta}_{g_n}$, where $\eo^*[\cdot] \equiv \eo[\cdot|\dd_n]$ denotes the bootstrap expectation conditional on the data $\dd_n \equiv \{(X_i, Y_i)\}_{i=1}^n$.
In other words,   $\hat{\ga}_{g_n}^{-1}$ does not generally match the inverse of $\hat{\ga}_n \equiv n^{-1}\sum_{i=1}^n (X_i - \bar{X})^{\otimes 2} $ for any finite truncation $g_n$.  However, by starting from an 
estimator $\hat{\beta}_{g_n} \equiv \hat{\ga}_{g_n}^{-1} \hat{\Dt}_n$, 
we may adjust a bootstrap-level estimating function to be   
\begin{align*}
	\Psi_i^*(\beta; \bar{X}, \bar{Y})
	& \equiv \check{\Psi}_i^*(\beta; \bar{X}, \bar{Y}) - \eo^*[\check{\Psi}_i^*(\hat{\beta}_{g_n}; \bar{X}, \bar{Y})]
	\\& = (X_i^* - \bar{X}) (Y_i^* - \bar{Y}) - (X_i^* - \bar{X})^{\otimes 2} \beta - \hat{U}_{n, g_n},
\end{align*}
by subtracting the bootstrap expectation
\begin{align} \label{eq_UnHatgn}
	\eo^*[\check{\Psi}_i^*(\hat{\beta}_{g_n}; \bar{X}, \bar{Y})] \equiv \hat{U}_{n,g_n} \equiv \frac{1}{n} \sum_{i=1}^n (X_i - \bar{X})(\hat{\e}_{i,g_n} - \bar{\hat{\e}}_{g_n}),
\end{align}
where the latter has a closed form expression based on the cross covariance between the regressors $\{X_i\}_{i=1}^n$ and the residuals $\{\hat{\e}_{i, g_n}\}_{i=1}^n$, $\hat{\e}_{i,g_n} \equiv Y_i - \langle \hat{\beta}_{g_n}, X_i \rangle$ arising  
from the estimator $\hat{\beta}_{g_n}$, with
$\bar{\hat{\e}}_{g_n} \equiv n^{-1} \sum_{i=1}^n \hat{\e}_{i,g_n}$ above. 
These corrected bootstrap estimating functions have bootstrap expectation of  
\begin{align*}
	\eo^*[\Psi_i^*(\beta; \bar{X}, \bar{Y})] = \hat{\Dt}_n  - \hat{\ga}_n \beta - \hat{U}_{n, g_n},
\end{align*}
which equals zero at $\beta = \hat{\beta}_{g_n}$ in the bootstrap world.  Consequently, $\hat{\beta}_{g_n}$ 
as the solution to $\eo^*[\Psi_i^*(\beta; \bar{X}, \bar{Y})]=0$ 
mimics true slope function $\beta = \ga^{-1} \Dt$  solving $\eo[\Psi_i(\beta; \mu_X, \mu_Y)] = 0$.

By replacing $\bar{X}$ and $\bar{Y}$ in $\Psi_i^*(\beta; \bar{X}, \bar{Y})$ with bootstrap data counterparts $\bar{X}^* \equiv n^{-1} \sum_{i=1}^n X_i^*$ and $\bar{Y}^* \equiv n^{-1} \sum_{i=1}^n Y_i^*$  (in analog to the original estimator $\hat{\beta}_{h_n}$ defined by using $\bar{X}$ and $\bar{Y}$ in place of $\mu_X$ and $\mu_Y$), 
a PB version $\hat{\beta}_{h_n}^*$ of the  FPCR estimator $\hat{\beta}_{h_n}$ is defined by the solution of the 
empirical bootstrap-data estimating equation  
\begin{align*}
	0=\frac{1}{n} \sum_{i=1}^n \Psi_i^*(\beta; \bar{X}^*, \bar{Y}^*) = \hat{\Dt}_n^* - \hat{\ga}_n^* \beta - \hat{U}_{n, g_n}, 
\end{align*}
upon regularization of  $(\hat{\ga}_n^*)^{-1}$, where 
$\hat{\ga}_n^* \equiv n^{-1} \sum_{i=1}^n (X_i^* - \bar{X}^*)^{\otimes 2}$ and $\hat{\Dt}_n^* \equiv n^{-1} \sum_{i=1}^n (Y_i^* - \bar{Y}^*) (X_i^* - \bar{X}^*)$ are averages from the bootstrap sample.  
Hence, the  PB re-creation of the FPCR estimator is then given by 
\begin{align}
	\hat{\beta}_{h_n}^* \equiv (\hat{\ga}_{h_n}^*)^{-1} (\hat{\Dt}_n^* - \hat{U}_{n, g_n}),
	\label{eq_fpcr_bts_bc}
\end{align}
{
	and the PB estimator of the intercept is then given by $\hat{\ap}_{h_n}^* \equiv \bar{Y}^* - \langle \hat{\beta}_{h_n}^*, \bar{X}^* \rangle$ in the same manner as for the original intercept estimator $\hat{\ap}_{h_n}$.
}
The construction in \eqref{eq_fpcr_bts_bc} matches how the original estimator $\hat{\beta}_{h_n}= \hat{\ga}_{h_n}^{-1}  \hat{\Dt}_n$ from \eqref{eq_fpcr} is the solution of $0=n^{-1}\sum_{i=1}^n \Psi_i(\beta; \bar{X}, \bar{Y})= \hat{\Dt}_n - \hat{\ga}_n \beta
$, based on \eqref{eqn:mest1}, upon similar regularization with truncation level $h_n$.  The combination  $(\hat{\beta}_{h_n}^*, \hat{\beta}_{g_n})$ in PB then serves to mimic  $(\hat{\beta}_{h_n},\beta)$ for inference about the FLRM.


\section{Distributional results under heteroscedasticity} \label{ssec_2_3}

\autoref{sssec_2_3_1}   describes a CLT for estimated conditional means $\hat{\mu}_{h_n}(X_0) \equiv \hat{\ap}_{h_n} + \langle \hat{\beta}_{h_n}, X_0 \rangle$ under the FLRM with  heteroscedasticity.  While  of potential interest in its own right, the CLT helps to develop 
the appropriate scaling needed for statistics  and to also frame some baseline assumptions that are useful for   
bootstrap.   \autoref{sssec_2_3_2} establishes the consistency of PB for distributional approximations.   For  comparison, \autoref{sssec_2_3_3} then provides a formal result to show that the naive implementation of bootstrap is generally invalid without restrictive conditions on truncation parameters.
{We close this section by discussing some theoretical aspects of our bootstrap results in \autoref{ssec_3_4}.}

\subsection{CLT for the conditional means under heteroscedasticity} \label{sssec_2_3_1}

Let $X_0$ denote a new regressor under the model, which is independent of data $\{(X_i, Y_i)\}_{i=1}^n$ and identically distributed as $X_1$. For an observed or given value of $X_0$ (i.e., conditional on $X_0$), 
{
we consider the sampling distribution of the difference 
\begin{align} 
	\label{eq_Sn}
	\sqrt{n \over s_{h_n}(X_0)} \{\hat{\mu}_{h_n}(X_0) - \mu(X_0)\},
\end{align}
between estimated  $\hat{\mu}_{h_n}(X_0) \equiv \hat{\ap}_{h_n} + \langle \hat{\beta}_{h_n},X_0\rangle$ and true $\mu(X_0) \equiv \ap + \langle \beta, X_0 \rangle$ conditional means.  
Above $s_{h_n}(X_0)$ denotes a  scaling factor, based on  $\ga_{h_n}^{-1} \equiv \sum_{j=1}^{h_n} \g_j^{-1} \pi_j$, which is defined as 
\begin{align} \label{eq_s}
	s_{h_n}(x) 
	& \equiv \langle \Ld \ga_{h_n}^{-1} (x-\eo[X]), \ga_{h_n}^{-1} (x-\eo[X]) \rangle  
	 = \|\Ld^{1/2} \ga_{h_n}^{-1} (x-\eo[X])\|^2
\end{align}
for $x \in \HH$, and involves the covariance operator $\Ld \equiv \eo[(X\e)^{\otimes 2}]$ of $X\e$, where $T^{1/2}$ denotes a self-adjoint square-root operator of a non-negative definite and bounded linear operator $T$ on $\HH$ such that $(T^{1/2})^2 = T^{1/2} T^{1/2} = T$.
A sample counterpart of \eqref{eq_s} is given as
\begin{align} \label{eq_sHat}
	\hat{s}_{h_n}(x) 
	& \equiv \|\hat{\Ld}_{n,k_n}^{1/2} \hat{\ga}_{h_n}^{-1} (x - \bar{X})\|^2, \quad x \in \HH,
\end{align}
where $\hat{\Ld}_{n,k_n} \equiv n^{-1} \sum_{i=1}^n \left( X_i \hat{\e}_{i,k_n} - n^{-1} \sum_{i=1}^n X_i \hat{\e}_{i,k_n} \right)^{\otimes 2}$ is an estimate of $\Ld$ based on residuals 
{$\hat{\e}_{i,k_n} \equiv Y_i - \hat{\mu}_{k_n} (X_i)$; }
for generality,  
here $k_n$ represents another tuning parameter used only to compute  residuals $\{\hat{\e}_{i,k_n}\}_{i=1}^n$ for estimated scaling 
$\hat{s}_{h_n}(x)$ in \eqref{eq_sHat}.  
}

Under either scaling factors $s_{h_n}(X_0)$ or $\hat{s}_{h_n}(X_0)$, we next show a CLT for the conditional mean $\mu(X_0) \equiv \ap + \langle \beta, X_0 \rangle$ in \autoref{thm_clt}, where the limiting distribution is   standard normal under the scaling.
For describing the CLT, some technical assumptions are listed.

\begin{enumerate}[(\text{A}1)] \itemsep 0cm
	\item \label{condA1_model}
	$\ker \ga = \{0\}$, where $\ker \ga \equiv \{x \in \HH: \ga x = 0 \}$ ; 
	
	\item \label{condA2_4thmoment}
	$\sup_{j \in \N} \g_j^{-2} \eo[ \langle X - \eo[X], \phi_j \rangle^4 ]<\infty$;
	
	\item \label{condA3_convexity}
	$\g_j $ is a convex function of $j \geq J$ (which implies that $  \g_j - \g_{j+1}$ is decreasing)  for some integer $J \geq 1$;
	
	\item \label{condA4_ev_rate}
	$\sup_{j \in \N} \g_j j \log j < \infty$;
	
	\item \label{condA5_eg_rate}
	$n^{-1} \sum_{j=1}^{h_n} \dt_j^{-2} \to 0$ as $n\to\infty$;
	
	\item \label{condA6_scale}
	$h_n s_{h_n}(X)^{-1} = O_\pr(1)$;
	
	\item \label{condA7_4thmoment_xe}
	there exists $\dt \in (0,2]$ such that $\sup_{j \in \N} \ld_j^{-(2+\dt)/2} \eo[|\langle X \e, \psi_j \rangle|]^{2+\dt} < \infty$,
	where $(\ld_j,\psi_j)$ is the $j$-th eigenvalue--eigenfunction pair of $\Ld$;
	
	\item \label{condA8_ga_ld} 
	$\sup_{j \in \N} \g_j^{-1} \|\Ld^{1/2} \phi_j\|^2 < \infty$.
	
\end{enumerate}

Condition \ref{condA1_model} is necessary for the model identifiability \citep{CFS99, CFS03, CMS07}. Conditions \ref{condA2_4thmoment} and \ref{condA7_4thmoment_xe} ensure that $X$ and $X\e$ respectively have finite fourth and $(2+\dt)$-th moments. 
Conditions \ref{condA3_convexity}-\ref{condA5_eg_rate} are   technical assumptions related to the decay behaviors of eigenvalues $\{\g_j\}$ and eigengaps $\{\dt_j\}$, where for \ref{condA4_ev_rate} we define $\dt_1 \equiv \g_1 - \g_2$ and $\dt_j \equiv \min \{ \g_j - \g_{j+1}, \g_{j-1} - \g_j \}$ for $j \geq 2$; such conditions are weak and are generally used to simplify proofs involving perturbation theory for functional data \citep{YDN23RB}.
Condition \ref{condA6_scale} provides a mild lower bound for scaling $s_{h_n}(X_0)$, where a similar assumption is needed in the homoscedastic setting  \cite{YDN23RB}.
Condition~\ref{condA8_ga_ld} is a technical condition that balances the eignedecay of $\ga$ and the decay rate of $\Ld$ in terms of $\{\phi_j\}_{j=1}^\infty$.
When Condition~\ref{condA2_4thmoment} holds, sufficient conditions for \ref{condA8_ga_ld} can also be developed by
assuming moment structures on the error and regressors; 
for example, Condition~\ref{condA8_ga_ld} follows 
if either $\eo[\e^4]<\infty$ 
or $\s^2(X) \equiv \eo[\e^2|X] = \sum_{j=1}^\infty \rho_j^2 \langle X, \phi_j \rangle^2$ for some $\{\rho_j\}_{j=1}^\infty$ such that $\sum_{j=1}^\infty \g_j \rho_j^2 < \infty$.
The statement of the CLT also involves 
the following condition,\\

\noindent Condition {$B(u)$}:
$\sup_{j \in \N} j^{-1} m(j,u) \langle \beta, \phi_j \rangle^2 < \infty$,\\

\noindent depending on a generic constant $u>0$ and function $m(j,u)$ of integer $j \geq 1$ defined as
\begin{align} \label{eq_bias_mju}
	m(j,u) = \max \left\{ j^u, \sum_{l=1}^j \dt_l^{-2} \right\}.
\end{align}

\noindent Condition $B(u)$ is generally mild and helps to remove bias in the limiting distribution of the statistics from \eqref{eq_Sn} by balancing the decay rates of eigenvalues and the Fourier coefficients of the slope function $\beta$, as described further in \autoref{rem_bias}.

A CLT for the conditional mean in FLRMs under heteroscedasticity is a new development in the FLRM literature, as given in the following theorem.

\begin{thm} \label{thm_clt}
	Suppose that Conditions \ref{condA1_model}-\ref{condA7_4thmoment_xe} hold along with $h_n^{-1} + n^{-1/2} h_n^{7/2} (\log h_n)^3 \to 0$ as $n\to\infty$. 
	We further suppose $n = O(m(h_n, u))$ along with Condition $B(u)$ for some $u>7$.
	Then, as $n\to\infty$, 
	\begin{enumerate}[(i)]
		\item 
		{
		\begin{align*}
			\sup_{y \in \R} \left| \pr \left( \sqrt{n \over s_{h_n}(X_0)} 
			\{\hat{\mu}_{h_n}(X_0) - \mu(X_0)\}
			\leq y \Big| X_0 \right) - \Phi(y) \right| \xrightarrow{\pr} 0,
		\end{align*}}
		where $\Phi$ denotes the standard normal distribution function.
		
		\item Additionally, if $\|\hat{\beta}_{k_n} - \beta\| \xrightarrow{\pr} 0$ and Condition \ref{condA8_ga_ld} hold, then  
		$\hat{s}_{h_n}(X_0)$ and $s_{h_n}(X_0)$ are asymptotically equivalent in that, for any $\eta>0$, 
		\begin{align*}
			\pr \left( \left| {\hat{s}_{h_n}(X_0) \over s_{h_n}(X_0)} - 1 \right| > \eta \Big| X_0 \right) \xrightarrow{\pr} 0,
		\end{align*}
		and  		 (i) also holds upon replacing $s_{h_n}(X_0)$ by the sample version $\hat{s}_{h_n}(X_0)$.
	\end{enumerate}
	
\end{thm}

\autoref{thm_clt}   generalizes the CLT for projections in FLRMs \citep{CMS07,YDN23RB} from the homoscedastic case  to broader heteroscedastic cases.   
When the errors are homoscedastic, i.e., $\eo[\e^2|X] \equiv \s_\e^2 \in (0, \infty)$, then the covariance operator of $\epsilon X$ becomes $\Ld = \s_\e^2 \ga$ and the scaling in \eqref{eq_sHat}  reduces to $s_{h_n}(X_0) = \s_\e^2 t_{h_n}(X_0)$, where  $t_{h_n}(x) \equiv \|\ga_{h_n}^{-1/2} x\|^2$ for $x \in \HH$. 
In this case,  \autoref{thm_clt} matches the CLT under homoscedasticity  \citep{CMS07, YDN23RB}.

From \autoref{thm_clt}, 
estimated conditional mean $\hat{\mu}_{h_n}(X_0) \equiv \hat{\ap}_{h_n} + \langle \hat{\beta}_{h_n},X_0\rangle$  with data-based scaling $\hat{s}_{h_n}(X_0)$ are asymptotically pivotal and, hence, a normal approximation may be applied to calibrate inference about true conditional mean $\mu(X_0) \equiv \ap + \langle \beta, X_0 \rangle$.
However,    resampling becomes useful for improving  distributional approximations in FLRMs under heteroscedasticity,  due to the complicated impacts of truncation $h_n$ in finite samples. The next section establishes the validity of   the proposed PB method.

\begin{rem}
	A sufficient condition for the consistency of $\hat{\beta}_{k_n}$ for $\beta$ in \autoref{thm_clt}(ii) is that $k_n^{-1} + n^{-1/2}k_n^2 \log k_n \to 0$ as $n\to\infty$; see also Theorem~S1, 
	\cite{supp}. 
	Theorem~\ref{thm_clt}(ii) may be further generalized by replacing $\hat{\beta}_{k_n}$   in the estimated scaling $\hat{s}_{h_n}(X_0)$ with a general consistent estimator of $\beta$.
\end{rem}


\begin{rem}
	Under   conditions on the error structure, 
	the rate on the truncation level $h_n$ can be weakened 
	to a lesser rate sufficient for obtaining a CLT under homoscedasticity.
	For instance, the rate $h_n^{-1}+ n^{-1/2} h_n^{5/2} (\log h_n)^2 \to 0$ is sufficient for \autoref{thm_clt}
	if either $\eo[\e^4]<\infty$ or $\eo[\e^2|X] = \sum_{j=1}^\infty \rho_j^2 \langle X, \phi_j \rangle^2$ for some $\{\rho_j\}_{j=1}^\infty$ with $\sum_{j=1}^\infty \g_j \rho_j^2<\infty$.
	This  rate matches ones assumed for the CLT  under homoscedasticity provided in \cite{CMS07,YDN23RB}.
	See also Remark~S1 in \cite{supp}. 
\end{rem}

\begin{rem} \label{rem_bias}
	In \autoref{thm_clt} and \autoref{thm_pb} to follow, the Conditions $n = O(m(h_n, u))$ (or $n = O(m(g_n, u))$) and $B(u)$ are necessary only for removing bias in limit distribution of 
	$\sqrt{n/s_{h_n}(X_0)} \{ \hat{\mu}_{h_n}(X_0) - \mu(X_0) \}$
	due to truncation $h_n$;  that is, without   these conditions, the asymptotic results would hold upon replacing $\mu(X_0) \equiv \ap + \langle  \beta, X_0 \rangle$ with a biased centering $\mu_{h_n}(X_0) \equiv \ap+ \langle \Pi_{h_n} \beta, X_0 \rangle$,
	where $\Pi_{h_n}$ denotes the projection on the first $h_n$ eigenfunctions $\{\phi_j\}_{j=1}^{h_n}$ of $\ga$ 
	and $\Pi_{h_n} \beta \equiv \sum_{j=1}^{h_n} \langle \beta, \phi_j \rangle \phi_j$ is a truncated version of the slope $\beta \equiv \sum_{j=1}^\infty \langle \beta, \phi_j \rangle \phi_j$.
	Such conditions are common for balancing the decay rates of eigenvalues and the Fourier coefficients of the slope function in the removal of bias $\langle (\Pi_{h_n} - I) \beta, X_0 \rangle$   from truncation;
	see \cite{CMS07, YDN23RB} for further discussion.

\end{rem}
{
\begin{rem} \label{remPoly}
	
	Related to the previous remark, 
	if we additionally assume polynomial decay rates as $\dt_j \equiv \g_j - \g_{j+1} \asymp j^{-a}$ and $|\langle \beta, \phi_j \rangle| \asymp j^{-b}$ for $a \in (2,\infty)$ and $b \in (5/2,\infty)$
	(though Condition~$B(u)$ with $u>7$ implies $\langle \beta, \phi_j \rangle = o(j^{-3})$),
	our balancing conditions along with rate of the tuning parameters can be reduced to the condition $h_n \asymp n^{1/v}$ for some $v \in (\max\{7, 2a+1\}, a+2b-1)$, which depends on the decay rates $a$ and $b$;
	here, for sequences $\{r_{1n}\}$ and $\{r_{2n}\}$ of positive real numbers, we use $r_{1n} \asymp r_{2n}$ to denote that ratios $r_{1n}/r_{2n}$ and $r_{2n}/r_{1n}$ are bounded away from zero.  Such specific decay rates of $\dt_j$ (or $\g_j$) and $\langle \beta, \phi_j \rangle$ are allowed and may provide a more familiar understanding of the growth rate conditions on truncation levels. 
	See \cite{CH06, HH07, lei14, LL22, YDN23RB, ZYZ23} for similar growth rates of truncation parameters.
\end{rem}
}

\subsection{Consistency of the paired bootstrap (PB)} \label{sssec_2_3_2}

Based on the CLT for   conditional means in \eqref{eq_SnHat}, we next consider PB approximations for the distribution of the studentized-type quantity
\begin{align} \label{eq_SnHat}
	T_n(X_0) \equiv \sqrt{n \over \hat{s}_{h_n}(X_0)} 
	{\{\hat{\mu}_{h_n}(X_0) - \mu(X_0)\}}
\end{align}
conditional on a given regressor $X_0$ with estimated scaling $\hat{s}_n(X_0)$  from \eqref{eq_sHat}. { 
A studentized bootstrap counterpart of \eqref{eq_SnHat}, with estimated bootstrap scaling $\hat{s}_{h_n}^*(X_0)$, is given as
\begin{align} \label{eq_SnHatStar_std}
	T_{n, \hat{s}^*}^*(X_0) \equiv \sqrt{n \over \hat{s}_{h_n}^*(X_0)} 
	{\{\hat{\mu}_{h_n}^*(X_0) - \hat{\mu}_{g_n}(X_0)\},}
\end{align} 
with the same fixed $X_0$, where $(\hat{\beta}_{h_n}^*,\hat{\beta}_{g_n})$  denote the bootstrap analogs \eqref{eq_fpcr_bts_bc} of the FPCR estimator $\hat{\beta}_{h_n}$ and true slope $\beta$.
Here, for $x \in \HH$, $\hat{\mu}_{h_n}^*(x) \equiv \hat{\ap}_{h_n}^* + \langle \hat{\beta}_{h_n}^*, x \rangle$ denotes the bootstrap estimated mean response,
where $\hat{\ap}_{h_n}^* \equiv \bar{Y}^* - \langle \hat{\beta}_{h_n}^*, \bar{X}^* \rangle$ with $\bar{X}^* \equiv n^{-1} \sum_{i=1}^n X_i^*$ and $\bar{Y}^* \equiv n^{-1} \sum_{i=1}^n Y_i$.
To define the bootstrap scaling $\hat{s}_{h_n}^*$ in \eqref{eq_SnHatStar_std}, recall that construction of  $\hat{s}_{h_n}$ in \eqref{eq_sHat} involves residuals from a FPCR estimator       
$\hat{\beta}_{k_n}$ with a generic bandwidth $k_n$.
A bootstrap version of scaling factor is then defined, in analog to   \eqref{eq_sHat}, as
\begin{align} \label{eq_sHatStar}
	\hat{s}_{h_n}^*(x)
	& \equiv \| (\hat{\Ld}_{n,k_n,g_n}^*)^{1/2} (\hat{\ga}_{h_n}^*)^{-1} (x-\bar{X}^*) \|^2, \quad x \in \HH,
\end{align}
where $\hat{\Ld}_{n,k_n,g_n}^* \equiv n^{-1} \sum_{i=1}^n ( X_i^* \hat{\e}_{i,k_n}^* - n^{-1} \sum_{i=1}^n X_i^* \hat{\e}_{i,k_n}^* )^{\otimes 2}$ is a bootstrap estimator of the covariance $\Ld$ based on  bootstrap residuals $\hat{\e}_{i,k_n}^* \equiv Y_i^* - \hat{\mu}_{k_n}^*(X_i^*) = Y_i^* - \hat{\ap}_{k_n}^* - \langle \hat{\beta}_{k_n}^*, X_i^* \rangle$ from bootstrap   estimators  $\hat{\ap}_{k_n} \equiv \bar{Y}^* - \langle \hat{\beta}_{k_n}^*, \bar{X}^* \rangle$ and $\hat{\beta}_{k_n}^* \equiv (\hat{\ga}_{k_n}^*)^{-1} (\hat{\Dt}_n^* - \hat{U}_{n, g_n})$; the latter is akin to  \eqref{eq_fpcr_bts_bc} with tuning parameter $k_n$.

Using the scaling factor $\hat{s}_n(X_0)$ in place of estimated scaling $\hat{s}^*_n(X_0)$ within bootstrap, another bootstrap counterpart of \eqref{eq_SnHat} can also be given as
\begin{align} 
\label{eq_SnHatStar}
	T_{n,\hat{s}}^*(X_0) \equiv \sqrt{n \over \hat{s}_{h_n}(X_0)} 
	{\{\hat{\mu}_{h_n}^*(X_0) - \hat{\mu}_{g_n}(X_0)\}.}
\end{align}
Due to the shared scaling in \eqref{eq_SnHat} and \eqref{eq_SnHatStar}, this bootstrap version essentially  approximates $ \hat{\mu}_{h_n}(X_0) - \mu(X_0)$ with $\hat{\mu}_{h_n}^*(X_0) - \hat{\mu}_{g_n}(X_0)$.

}

{
In  formal results next, the validity of the bootstrap in FLRMs depends on the asymptotic relationship between truncation levels $h_n$ and $g_n$ through   the following limiting ratio:
\begin{align} \label{eq_ratio_hg}
	\tau \equiv \lim_{n\to\infty} h_n/g_n .
\end{align}
We   show that our modified paired bootstrap is generally flexible
in the sense that  bootstrap consistency holds even when $h_n$ is asymptotically bigger than $g_n$, i.e., $\tau \geq 1$ (cf.~\autoref{thm_pb});
in contrast, the naive paired bootstrap is more restricted, as described in \autoref{sssec_2_3_3}.
}

\autoref{thm_pb} establishes the consistency of the PB method for the sampling distribution of the studentized conditional mean estimator   in \eqref{eq_SnHat} under heteroscedasticity. 
Let $\pr^* \equiv \pr(\cdot|\dd_n)$ denotes the bootstrap probability conditional on the sample $\dd_n \equiv \{(X_i, Y_i)\}_{i=1}^n$.   

\begin{thm} \label{thm_pb}
	Suppose that Conditions \ref{condA1_model}-\ref{condA8_ga_ld} hold,
	that $k_n^{-1} + n^{-1/2} k_n^2 \log k_n \to 0$  as $n\to\infty$,
	and that $\eo[\|X\|^{4+2\dt}] <\infty$ and $n^{-\dt/2} h_n^{\dt/2} \sum_{j=1}^{h_n} \ld_j^{-(2+\dt)/2} = O(1)$ as $n\to\infty$ for $\dt \in (0,2]$ in Condition~\ref{condA7_4thmoment_xe}.
	Along with Condition $B(u)$ for some $u>7$, we further suppose that {the limiting ratio $\tau$ from \eqref{eq_ratio_hg} is not less than 1,} $g_n^{-1} + n^{-1/2} h_n^{7/2} (\log h_n)^3 \to 0$, and $n = O(m(h_n, u))$.
	Then, as $n\to\infty$, the paired bootstrap (PB) is valid for the distribution of the studentized statistic $T_n(X_0)$ in \eqref{eq_SnHat}:
	\begin{align*}
		\sup_{y \in \R} \left| 
		\pr^* \left( T_n^*(X_0) \leq y | X_0 \right)
		- \pr \left( T_n(X_0) \leq y |  X_0\right)
		\right| \xrightarrow{\pr} 0,
	\end{align*}
	where $T^*_n(X_0)$  denotes either  $T_{n,\hat{s}^*}^*(X_0)$ from \eqref{eq_SnHatStar_std} or   $T_{n,\hat{s}}^*(X_0)$  from \eqref{eq_SnHatStar}.
\end{thm}

\autoref{thm_pb} conditions for the PB are similar to those for the CLT   from  \autoref{thm_clt}, 
though additional mild assumptions (i.e., $\tau \equiv \lim_{n\to\infty}h_n/g_n \geq 1$) appear to connect the second truncation $g_n$  in PB to the original data truncation $h_n$.     
Namely, the truncation level $g_n$  for defining the bootstrap rendition $\hat{\beta}_{g_n}$ of the true parameter $\beta$  may differ from the truncation  $h_n$ used in the original FPCR estimator $\hat{\beta}_{h_n}$, 
though $g_n$ may not be larger than $h_n$ asymptotically (see also \autoref{prop_pb_both_fail}).  
This coordination of truncation levels is generally required for the bootstrap to be asymptotically correct, which allows the bootstrap to control the bias type described in \autoref{rem_bias}.
In practice, we recommend choosing a slightly smaller $g_n$ than $h_n$.
In particular, we give  a rule of thumb for selecting $h_n$ and $g_n$  in \autoref{sec4}, which performs well as illustrated numerically.

{

\begin{eg} \label{eg1}
	As a data-generation for illustrating Theorem~\ref{thm_clt}-\ref{thm_pb},
	we consider  $\HH = L^2([0,1])$, 
	the space of square integrable functions defined on the interval $[0,1]$, 
	and suppose the eigenfunctions $\{\phi_j\}$ are Fourier basis functions.
	We further suppose the following distributional conditions (a)-(b) on $X$ and $\e$:  
	\begin{enumerate}[(a)]
		
		\item 
		the conditional variance of the error $\e$ given the regressor $X$ is $\s^2(X) \equiv \eo[\e^2|X] = \sum_{j=1}^\infty \rho_j^2 \langle X, \phi_j \rangle^2$ for some $\{\rho_j\}_{j=1}^\infty$ with $\sum_{j=1}^\infty \g_j \rho_j^2 < \infty$;
		
		\item 
		$X$ has functional principal component (FPC) scores as $\g_j^{-1/2} \langle X, \phi_j \rangle = \xi W_j$ for $j\geq 1$, where $\{W_j\}$ denote   iid standard normal variables and, independently,  
		$\xi$ is a general random variable.
		
	\end{enumerate}
As in numerical studies later (Section~\ref{ssec_4_1}), suppose further  that the error distribution is a centered $\chi^2$ distribution as $\e|X \sim \chi^2(\nu(X)) - \nu(X)$ with $\nu(X) = \|X\|^2/2$ (i.e., the case with $\rho_l=1$ for all $l$) and $\eo[\xi^{10}]<\infty$.  
	Under polynomial decay rates $\dt_j \equiv \g_j- \g_{j+1} \asymp j^{-a}$ and $|\langle \beta, \phi_j \rangle| \asymp j^b$ with $a>2$ and $b>1$,
	let the truncation level $h_n$ grow  as $h_n \asymp n^{1/v}$ where $7 \wedge (2a+1) < v < a+2b-1$ (cf.~\autoref{remPoly}), and $\tau \geq 1$ from \eqref{eq_ratio_hg}.  Then, all conditions of Theorem~\ref{thm_clt}-\ref{thm_pb} may be verified to hold;
	see \autoref{app3} for the proof.
\end{eg}
}

{
  \autoref{thm_pb} can also be  slightly modified  for inference about the projection $\langle \beta, X_0 - \eo[X] \rangle$ of the slope function $\beta$ onto a centered new predictor $X_0 - \eo[X]$.   
This
can be useful for isolating
  the  effect of the slope $\beta$ on the conditional mean $\mu(X_0) = \ap + \langle \beta, X_0 \rangle = \E[Y] + \langle \beta, X_0 - \eo[X] \rangle$, where   estimation of the intercept $\ap$ is not involved.
The projection statistics and bootstrap approximations in \autoref{cor_proj} to follow can additionally play a role in developing   hypothesis tests about multiple   projections, as described in \autoref{sec3}.  
Namely, consider approximating the distribution of the studentized projection statistic
\begin{align} \label{eq_stat_proj}
    T_{n, proj}(X_0) & \equiv \sqrt{n \over \hat{s}_{h_n}(X_0)} (\langle \hat{\beta}_{h_n}, X_0 - \bar{X} \rangle - \langle \beta, X_0 - \eo[X] \rangle)
\end{align}
with a bootstrap counterpart as
\begin{align} \label{eq_stat_proj_bts}
    T_{n, \hat{s}^*,proj}^*(X_0) & \equiv \sqrt{n \over \hat{s}_{h_n}^*(X_0)} (\langle \hat{\beta}_{h_n}^*, X_0 - \bar{X}^* \rangle - \langle \hat{\beta}_{g_n}, X_0 - \bar{X} \rangle)
\end{align}
or as
\begin{align} \label{eq_stat_proj_bts1}
    T_{n, \hat{s},proj}^*(X_0) & \equiv \sqrt{n \over \hat{s}_{h_n}(X_0)} (\langle \hat{\beta}_{h_n}^*, X_0 - \bar{X}^* \rangle - \langle \hat{\beta}_{g_n}, X_0 - \bar{X} \rangle),
\end{align}
where \eqref{eq_stat_proj_bts}-\eqref{eq_stat_proj_bts1} 
are versions of \eqref{eq_SnHatStar_std} and \eqref{eq_SnHatStar} for centered projections.
\autoref{cor_proj} then justifies the bootstrap  with  projection-based quantities \eqref{eq_stat_proj}-\eqref{eq_stat_proj_bts1}.

\begin{cor} \label{cor_proj}
    Under the  assumptions of \autoref{thm_pb},  as $n\to\infty$, the PB is   valid for approximating the distribution of the projection statistic $T_{n,proj}(X_0)$    in \eqref{eq_stat_proj}:
    \begin{align*}
        \sup_{y \in \R} |\pr^*(T_{n, proj}^*(X_0) \leq y|X_0) - \pr(T_{n, proj}(X_0) \leq y|X_0)| \xrightarrow{\pr} 0,
    \end{align*}
   where $T^*_{n,proj}(X_0)$  denotes either the bootstrap rendition from \eqref{eq_stat_proj_bts} or from   \eqref{eq_stat_proj_bts1}.
\end{cor}
When regressor and response variables have   mean zero ($\eo[X] =0=\eo[Y]$), as assumed in some FLRM works \cite{CH06, CMS07, GM11, HH07,  lei14, YDN23RB, ZYZ23}, 
then  \autoref{cor_proj} remains valid 
by simply dropping  regressor  means $\eo[X], \bar{X}, \bar{X}^*$ in \eqref{eq_stat_proj}-\eqref{eq_stat_proj_bts1}.



}

\subsection{Limitations of naive bootstrap} \label{sssec_2_3_3}

As described in \autoref{ssec_2_2}, a naive bootstrap formulation $\hat{\beta}_{h_n,naive}^* \equiv (\hat{\ga}_{h_n}^*)^{-1} \hat{\Dt}_n^*$ of the FPCR estimator  will not be generally be valid for approximating the distribution a conditional mean statistic  
$T_n \equiv  \sqrt{ n/ \hat{s}_{h_n}(X_0)} \{\hat{\mu}_{h_n}(X_0) - \mu(X_0)\}$ 
in \eqref{eq_SnHat} {\it unless} bootstrap centering parameter $\beta^*\equiv \hat{\beta}_{g_n}$ is narrowly chosen.  
{That is, if the limiting ratio $\tau $ in \eqref{eq_ratio_hg} is bigger than 1,
the naive bootstrap can fail,
whereas the PB method of \autoref{sssec_2_3_2}  is consistent.}
 As a formal illustration,  we     consider   { a bootstrap quantity 
 \begin{align} \label{eq_SnStarCheck}
	T_{n,naive}^*(X_0) \equiv \sqrt{n \over \hat{s}_{h_n,naive}^*(X_0)} 
	\{\hat{\mu}_{h_n, naive}^*(X_0) - \hat{\mu}_{g_n}(X_0)\} 
\end{align}
that differs from a valid bootstrap version $T_{n,\hat{s}^*}^*(X_0)$  in   \eqref{eq_SnHatStar_std}  by using
$\hat{\mu}_{h_n, naive}^*(X_0) \equiv \bar{Y}^* + \langle \hat{\beta}_{h_n,naive}^*, X_0 - \bar{X}^*\rangle$ with the naive bootstrap estimator  $\hat{\beta}_{h_n,naive}^*$ in  place of the proposed PB $\hat{\beta}_{h_n}^*$; the bootstrap scaling $\hat{s}_{h_n,naive}^*(X_0)$ in \eqref{eq_SnStarCheck} uses  $\hat{\beta}_{k_n,  naive}^*=(\hat{\ga}_{k_n}^*)^{-1} \hat{\Dt}_n^*$ when computing the residuals for \eqref{eq_sHatStar}. }
\autoref{prop_pb_naive_fail} provides a general illustrative data example where the naive bootstrap method provably fails, which stands in contrast to the consistency of the modified PB   from    \autoref{thm_pb}.
In the following, distributional convergence is  in
the space of  real-valued functions on
$[-\infty,\infty]$ that are right continuous with left limits  (cf.~\cite{bill99}).

\begin{prop} \label{prop_pb_naive_fail}
	Suppose conditions (a)-(b) {in \autoref{eg1} hold with $\eo[\xi^8]<\infty$};  that \autoref{thm_pb} assumptions  hold with $\dt = 2$; 
 that $n^{-1/2} h_n^{9/2} (\log h_n)^6 \to 0$; and  {that the limiting ratio $\tau$ from \eqref{eq_ratio_hg} is larger than 1}.  
%
%
%
	Then,  
	as $n\to \infty$,
	\begin{align*}
		\pr^* ( T_{n,naive}^*(X_0) \leq y | X_0 ) - \pr(T_{n}(X_0) \leq y|   X_0)   \xrightarrow{d} \Phi\Big(y+\s(\tau) Z \Big) - \Phi(y), \; y\in\mathbb{R},
	\end{align*}
	  where $Z$ denotes a standard normal variable with distribution function $\Phi$ and $\s(\tau)>0$ denotes a constant (cf.~(\ref{eq_consbias_limvar})).   Thus, the naive bootstrap is  inconsistent.
\end{prop}

	

\begin{rem} \label{remnew}
	The naive bootstrap \eqref{eq_SnStarCheck} in \autoref{prop_pb_naive_fail}   can be shown to be  valid upon restricting  {$\tau$ from \eqref{eq_ratio_hg} to equal   1},
	which in case $\s(\tau)=0$ holds (cf.~(\ref{eq_consbias_limvar})) and the limit becomes zero. Essentially, {for the naive   bootstrap to work,} the bootstrap centering $\beta^*\equiv\hat{\beta}_{g_n}$ must be confined to the original FRCR estimator $\hat{\beta}_{h_n}$ (i.e., $g_n=h_n$).  Further, neither the proposed or naive PB approach is generally valid 
	if {$\tau<1$} (cf.~\autoref{prop_pb_both_fail}).
\end{rem}

\begin{figure}[!b]
	\centering
	\includegraphics[width=0.79\linewidth]{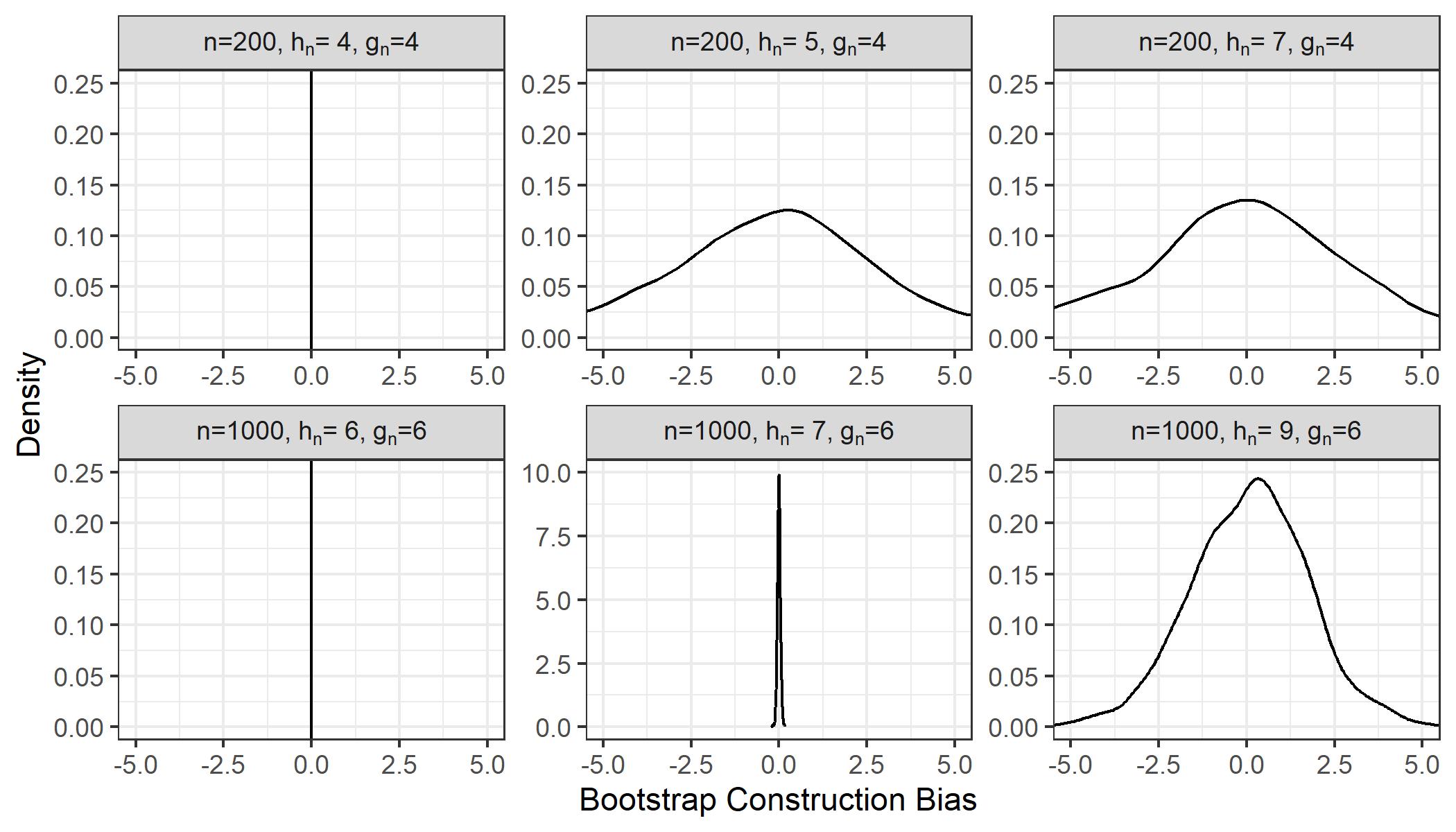}
	\caption{
		Kernel density estimates of a construction bias  $B_n$ in  the naive bootstrap. Plots have a common $x$-axis, and  the bias is zero   when $g_n=h_n$.
	}
	\label{fig_PBfail_cons}
\end{figure}

A take-away from \autoref{prop_pb_naive_fail} is that the naive bootstrap can fail with simple regressor structures, such as Gaussian $X$ (i.e., $\xi=1$ above), though Condition~(b) of \autoref{eg1} serves to accommodate a larger class of regressor distributions with potential dependence among FPCs.      That is, the  distance 
between the true distribution   of the target statistic $T_n(X_0)$ and that of the naive bootstrap approximation  $T_{n,naive}^*(X_0)$ can have a random limit and may not converge to zero at any point  on the real line. 
This aspect arises due to an extra construction bias created 
in the naive bootstrap definition of $\hat{\beta}_{h_n,naive}^*$ under heteroscedasticity,  which may be explained as follows.  
From the proof of \autoref{prop_pb_naive_fail},  quantiles from the naive bootstrap approximation  $T_{n,naive}^*(X_0)$ in \eqref{eq_SnStarCheck} are shifted  from those of a valid bootstrap approximation $T_{n,\hat{s}}^*(X_0)$ in \eqref{eq_SnHatStar}
by a   random  contribution, say $B_n$, that depends on the original data but not the bootstrap sample;   
in large samples, this bias amount $B_n\approx T_{n,naive}^*(X_0)-T_{n,\hat{s}}^*(X_0)$ acts like a draw  from a mean-zero normal distribution having a variance
\begin{align} \label{eq_consbias_limvar}
	\s^2(\tau) \equiv (1-\tau^{-1}) \left( \|\ga^{1/2}\beta\|^2 \Big/ \sum_{j=1}^\infty \g_j \rho_j^2 + 1 \right)
\end{align}
{that is non-zero when the   ratio $\tau$ from \eqref{eq_ratio_hg} exceeds 1}.   
This bias behavior can also be observed practically.
\autoref{fig_PBfail_cons} contains a numerical illustration 
based on 1000 experiments generated from a FLRM with regressor $X$ and error $\e$ as described in \autoref{prop_pb_naive_fail}.
We examine the resulting distribution of the construction bias $B_n$   in the naive approach when $h_n/g_n > 1$.
\autoref{fig_PBfail_cons} shows the distribution of this   term $B_n$ is quite different from zero as the ratio $h_n/g_n$ increases, even when $h_n = g_n + 1$, 
so that this bias $B_n$ is non-ignorable.
The latter observation matches the theoretical result in \eqref{eq_consbias_limvar}, underling    \autoref{prop_pb_naive_fail}.

For clarity, both naive and modified PB may fail
{if the limiting ratio $\tau$ from \eqref{eq_ratio_hg} is less than 1}
due to a different source of bias (i.e.,  apart from the   construction of the bootstrap estimator $\hat{\beta}_{h_n}^*$), which relates to    centering  in the CLT (cf.~\autoref{rem_bias}).  
This bias does not vanish if $h_n/g_n<1$, which arises because any estimator $\hat{\beta}_{g_n}$   playing the role of the true slope $\beta$ in the bootstrap world cannot capture the infinite dimensionality of $\beta$.
This failure is illustrated in \autoref{prop_pb_both_fail}, 
with details in the supplement \cite{supp}.

\begin{prop} \label{prop_pb_both_fail}
	Suppose that the assumptions of \autoref{thm_pb} with $\dt = 2$ hold and that $\tau \in (0,1)$. 
	We further suppose Conditions (a) and (b) in \autoref{eg1}.
	Then, as $n\to\infty$, both naive $T_{n,naive}^*(X_0)$ and modified $T_n^*(X_0)$ bootstrap renditions of  $T_n(X_0)$ fail to provide  asymptotically correct distributional approximations. 
	Namely, the  convergence   in \autoref{prop_pb_naive_fail}  holds for $T_{n,naive}^*(X_0)$  upon redefining     $\s(\tau)$ there as $\sqrt{\tau^{-1}-1}>0$, and this  same result holds also for $T_n^*(X_0)$.
\end{prop}

To summarize, 
{the asymptotic ratio $\tau \in (0, \infty)$ from \eqref{eq_ratio_hg}} plays a significant role in   PB methods.
The naive PB requires asymptotic equivalence of $h_n$ and $g_n$ with $\tau=1$ and becomes invalid when $h_n$ is asymptotically larger than $g_n$ with $\tau>1$.  
In contrast, the modified PB enjoys additional flexibility in setting $h_n$ and $g_n$ with $\tau \geq 1$, which results in better numerical performance in \autoref{sec4}. Simulation evidence indicates that	our modified PB produces good and stable coverages over different $h_n \geq g_n$,
while   the naive PB is more sensitive to the ratio $h_n/g_n$ and tends to over-cover.
Both modified and naive PB methods can fail  if $h_n$ is asymptotically smaller than $g_n$ with $\tau<1$.

{
\subsection{Theoretical discussions and extensions} \label{ssec_3_4}

We provide some context  about theoretical differences   that distinguish our work from previous resampling with FLRMs  (cf.~\cite{YDN23RB}).
The majority of challenges owe to the heteroscedastic structure of the data. 
 To our knowledge, there has been no theoretical development of  asymptotic inference  in heteroscedastic FLRMs, while existing work considers only homoscedastic errors.\citep{CMS07, GM11,  KH16a, KH16b, YDN23RB}.

Similar to  previous general works on FLRMs \cite{CH06,  CMS07, HH07, GM11}, our proposed PB and previous developments of residual bootstrap (RB) \cite{YDN23RB} consider random   regressors $\{X_i\}_{i=1}^n$
and account for this randomness in inference.    Under homoscedastic errors, RB can also provide inference about conditional sampling distributions given  regressors, which is true in finite-dimensional linear models \cite{free81} and FLRMs \cite{YDN23RB}.  However, a larger point is that the RB for FLRMs has no guarantees for  inference outside of homoscedastic errors, while the PB is generally valid under much more general error distributions, including heteroscedasticity.  Numerical studies in \autoref{sec4} demonstrate that the PB can perform similarly to RB in coverage accuracy under homoscedasticity, where RB is expected to have advantages, but can largely outperform RB when dependence exists between errors and regressors.

Compared to previous work with FLRMs, one large complication here has involved   identifying the appropriate scaling terms $s_{h_n}$ as in \eqref{eq_s} so that estimated conditional means/projections  have well-defined  limits     under heteroscedasticity (\autoref{thm_clt}).  Such scaling is more complicated compared to the homoscedastic case, and  
the consistency of   the sample scaling $\hat{s}_{h_n}$ (\autoref{thm_clt}) as well as the bootstrap counterpart  $\hat{s}_{h_n}^*$ (\autoref{thm_pb})  become non-trivial to establish (cf.~supplement \cite{supp}).  Further technical complications also arise in bootstrap-level versions of perturbation theory (cf.~\cite[][Section~S3.1]{supp}).
Namely, due to the resampling nature of the PB (i.e., resampling both regressors and responses),
standard results on perturbation theory (e.g.,~\cite[][Chapter~5]{HE15}) do not directly apply to the bootstrap quantities,   requiring a different treatment that does not arise in previous work 
on RB 
  \cite{YDN23RB}.  Another non-trivial aspect of PB for FLRMs under heteroscedasticity is the inconsistency of the naive PB due to a construction bias 
  (\autoref{prop_pb_naive_fail}); a similar  bias does not arise for PB in
  standard finite-dimensional regression  \cite{free81}   or for RB in    homoscedastic   FLRMs \cite{YDN23RB}

}

{
\begin{rem} \label{rem_nonfully}
	
	In our development, we essentially consider only the case where functional data can be fully observed.
	When functional data are not fully observed, 
	the bootstrap theory would require
   modification   due to several complicating factors.
	The covariance estimator $\hat{\ga}_n \equiv n^{-1} \sum_{i=1}^n (X_i - \bar{X})^{\otimes 2}$ may not applicable with irregular (and possibly sparse) time grid points, so that the FPCR estimator may not be directly used.
	Additionally,  our results rely heavily on a perturbation theory and  its bootstrap version  developed in the current paper,
 though these would require adjustments for noisy and irregularly observed functional data  (e.g., see \cite{ZYZ23} for a recent non-bootstrap version of perturbation theory for non-fully observed functional data).
	   Related to this last point, 
    while this work provides a first CLT for conditional means in FLRMs under heteroscedasicity,  
    an alternative  development for CLT may be needed for irregularly observed functional regressors as perturbation theory often plays a  key role with bias terms arising in CLTs (cf.~\cite{CMS07, YDN23RB}). 
	Bootstrap methods may then also require re-formulation for irregular/sparse functional data.   
	This is beyond our current scope, which we mention for future research.
\end{rem}
}

{
 
\begin{rem} \label{rem_noniid}
 The approximation results here can also extend to a more general new predictor $X_0$, 
    which need not be independent of the data or share the same distribution as an underlying regressor function $X$. 
    For instance, our results remain valid
   if $X_0$ follows a different distribution
  than $X$ but with matching mean and variance (or eigenvalues $\{ \mu_j \}$ of $\var[X_0]$ can match  those $\{\g_j\}$ of $\var[X]$ asymptotically).
    Another simple scenario involves $X_0$ as dependent on the data regressors $\{X_i\}_{i=1}^n$, such as $X_0 = X_1$ or an average from some subset of $\{X_i\}_{i=1}^n$.  See also \cite{YDN23RB} for a related discussion about RB.
\end{rem}

}


\section{Hypothesis tests for multiple conditional means} \label{sec3}

{
The testing of the association between the functional regressor $X$ and the scalar response $Y$ in FLRMs has drawn much recent attention in the literature.
    A global test of $\beta=0$ was first proposed in \cite{CFMS03} by assessing the covariance operator for   $\Dt = 0$,   
and   several  works have similarly considered various global tests     \citep{CGS04,GGMG12, HMV13, lei14, KSM16, SDH17, LL22}.  For instance, in \cite{KSM16, SDH17}, the global null $\beta = 0$ is approximated by testing $\langle \beta, \phi_1 \rangle = \cdots = \langle \beta, \phi_L \rangle = 0$  for an increasing integer $L$ depending on $n$, where $\{\phi_j\}_{j=1}^\infty$ denotes the set of eigenfunctions of $\ga \equiv \eo[(X-\eo[X])^{\otimes 2}]$; this approach     essentially involves a linear model with increasing number of scalar parameters, whereby the (sample) FPC scores are treated as the observed regressors  rather than  $\{X_i\}_{i=1}^n$. 
However, none of these previous works  applies to  testing  claims regarding the orthogonality of $\beta$ to general regressor functions observed from random regressors.  The proposed PB method from  \autoref{sec2}, though, can be adapted  to assess projections of the slope function $\beta$  onto subspaces spanned by general directions, as next explained.   

}

{
To frame the testing problem,  let $\xx_0  \equiv \{X_{0,l}\}_{l=1}^L$  denote an observed collection of random regressors   under consideration, for some   $L\geq 1$. 
We wish to test the null hypothesis about  constancy of the collective  conditional means as
\begin{align}
    H_0: \mu(X_{0,l}) = \eo[Y], ~\forall l =1, \dots, L \quad \mathrm{against} \quad H_1: H_0 \mathrm{~is~not~true} \label{eq_null}
\end{align}
The null hypothesis states that the linear effects of the slope function $\beta$   do not change across the  new regressors $\xx_0  \equiv \{X_{0,l}\}_{l=1}^L$ and, in ANOVA fashion, equal a common average $\mu(\eo[X])=\eo[Y]$ set by the global mean curve $\eo[X]$. The null hypothesis in \eqref{eq_null} can also be  expressed in terms of projections as 
\begin{align} \label{eq_null_proj}
    H_0: \langle \beta, X_{0,l} \rangle = \langle \beta, \eo[X] \rangle, ~ \forall l=1, \dots, L \quad \mathrm{against} \quad H_1: H_0 \mathrm{~is~not~true}.
\end{align}
In particular cases with zero mean regressors (i.e., $\eo[X] = 0$), the null hypothesis \eqref{eq_null_proj}
boils down to assessing
the orthogonality of the slope function $\beta$ to the linear subspace $\mathrm{span} (\xx_0) \subseteq \HH$ spanned by $\xx_0$, or 
$\Pi_{\xx_0}\beta = 0$ where $\Pi_{\xx_0} \beta$  denotes the projection   of $\beta$ onto $\mathrm{span} (\xx_0)$.
Due to the equivalence between \eqref{eq_null} and \eqref{eq_null_proj}, the test statistics based on the projection quantities \eqref{eq_stat_proj}-\eqref{eq_stat_proj_bts} can be applied to test   null hypothesis \eqref{eq_null}
without estimation of the intercept $\ap$.
}
As the bootstrap results in \autoref{ssec_2_3} apply for a given regressor $X_0$,  
a PB-based testing procedure  can be formulated to assess this type of hypothesis. 
An advantage is that this approach  provides a specific test of whether regression effects exist in any pre-defined directions,
while a global test about $\beta$ (e.g.,  \cite{KSM16, SDH17}) is not amenable to this purpose.
Additionally, previous works on hypothesis testing for FLRMs rely heavily on independent regressor-error assumptions, while our bootstrap-based testing procedure can address  such testing problems in FLRMs under conditionally dependent errors and heteroscedasticity.

{
To describe test statistics, write
\begin{equation}
	\label{eqn:Hostat}
	{T}_{n,l}^{\scriptscriptstyle H_0} \equiv \sqrt{ \frac{n} {\hat{s}_{h_n}(X_{0, l})}} 
	 \langle \hat{\beta}_{h_n}, X_{0, l} - \bar{X} \rangle,
	\quad l=1,\ldots,L,
\end{equation}
 to denote the studentized projection statistic from \eqref{eq_stat_proj}  upon substituting the hypothesized value $\langle \beta, X_{0,l} -\eo[X] \rangle=0$  by \eqref{eq_null_proj} for a new regressor $X_{0, l}$. }  
We may define test statistics by combining these  direction-based statistics   as     
\begin{align} \label{eq_testing_stat}
	W_{n,L^2} \equiv \sum_{l=1}^L \left[{T}_{n,l}^{\scriptscriptstyle H_0}\right]^2
	\quad \text{and} \quad
	W_{n, \max} \equiv \max_{1 \leq l \leq L} \left|{T}_{n,l}^{\scriptscriptstyle H_0}\right|, 
\end{align}
representing $L_2$- or $L_\infty$-type norms.  Large values of such statistics then provide evidence against  $H_0$.   These test statistics are well-defined with  non-degenerate limit distributions under the null hypothesis, though their limit laws are complicated, depending intricately on covariances between estimated projections.  Consequently, these limit distributions are impractical for direct use.  However, the sampling distributions of test statistics can be viably approximated with the proposed PB method and, in fact, there exist two ways of implementing the bootstrap here: by enforcing the null hypothesis at the bootstrap level or not.

If we do not enforce the null hypothesis in the bootstrap world, then we essentially adopt the same   PB procedure described in  \autoref{sssec_2_3_2}.  That is, we may formulate studentized bootstrap quantities for centered projections, similar to \eqref{eq_stat_proj_bts}, as 
\begin{align*}
	T_{n,l,\hat{s}^*}^* \equiv \sqrt{n \over \hat{s}_{h_n}^*(X_{0, l})} 
	{\{ \langle \hat{\beta}_{h_n}^*, X_{0, l} - \bar{X}^* \rangle - \langle \hat{\beta}_{g_n}, X_{0, l} -\bar{X} \rangle \}},
	\quad l=1,\dots,L
\end{align*}
based on the same bootstrap sample $\{(X_i^*,Y_i^*)\}_{i=1}^n$ and a common estimator $\hat{\beta}_{g_n}$ playing the bootstrap role of $\beta$.  The bootstrap test statistics are then   
\begin{align} \label{eq_testing_statHat}
	W_{n,L^2}^* \equiv \sum_{l=1}^L \left[ T_{n,l,\hat{s}^*}^*\right]^2 
	\quad \text{and} \quad 
	W_{n, \max}^*\equiv \max_{1 \leq l \leq L} \left| T_{n,l,\hat{s}^*}^*\right|, 
\end{align}

To enforce the null hypothesis in the bootstrap world, we  modify the PB procedure described in \autoref{ssec_2_2}, letting $\tilde{\beta}_{g_n} \equiv \hat{\beta}_{g_n} - \Pi_{\xx_0^c} \hat{\beta}_{g_n}$ rather than $\hat{\beta}_{g_n}$
denote the bootstrap analog of the slope $\beta$.  Here $\tilde{\beta}_{g_n}$ denotes a version of $\hat{\beta}_{g_n}$ after removing its projection $\Pi_{\xx_0^c} \hat{\beta}_{g_n}$ onto the subspace spanned by $\xx_0^c \equiv \{ X_{0,l} - \bar{X} \}_{l=1}^L $.  With this change,    
it holds that $\Pi_{\xx_0^c} \tilde{\beta}_{g_n} = 0$ and so $\tilde{\beta}_{g_n}$ mimics the  same property  $\langle \beta, X_{0,l} - \eo[X] \rangle = 0$, $l=1, \dots, L$  of the true parameter $\beta$ under $H_0$.  To formulate bootstrap data, we also 
write a response variable $\tilde{Y}_i \equiv Y_i - \langle \Pi_{\xx_0^c} \hat{\beta}_{g_n}, X_i \rangle$ after removing a projection contribution from $\Pi_{\xx_0} \hat{\beta}_{g_n}$ with respect to $X_i$.  A PB sample $\{(X_i^*, \tilde{Y}_i^*)\}_{i=1}^n$ is  defined by iid draws from the empirical distribution of $\{(X_i, \tilde{Y}_i)\}_{i=1}^n$, and 
the same development from \autoref{ssec_2_2} then applies with the change that $Y_i^*$, $\hat{\beta}_{g_n}$, $\bar{Y} $ there   become   $\tilde{Y}_{i}^*$, $\tilde{\beta}_{g_n}$, $\bar{\tilde{Y}} \equiv n^{-1} \sum_{i=1}^n \tilde{Y}_i = \bar{Y} - \langle \Pi_{\xx_0^c} \hat{\beta}_{g_n}, \bar{X} \rangle$. 
The bootstrap estimator then has a  
closed form as  
\[
\tilde{\beta}_{h_n}^* \equiv (\hat{\ga}_{h_n}^*)^{-1} (\tilde{\Dt}_n^* - \hat{U}_{n, g_n})
\]
in parallel to \eqref{eq_fpcr_bts_bc} with $\tilde{\Dt}_n^*
\equiv n^{-1} \sum_{i=1}^n (\tilde{Y}_i^* - \bar{\tilde{Y}}^*) (X_i^* - \bar{X}^*)$  in place of   $\hat{\Dt}_n^* \equiv n^{-1} \sum_{i=1}^n (Y_i^* - \bar{Y}^*) (X_i^* - \bar{X}^*)$.   
When enforcing the null hypothesis at the bootstrap level,
bootstrap versions  of  test statistics  in \eqref{eq_testing_stat} are then given by \begin{align} \label{eq_testing_statTilde}
	W_{n,L^2}^* \equiv \sum_{l=1}^L \left[{T}_{n,l}^{* {\scriptscriptstyle H_0}}\right]^2 
	\quad \text{and} \quad 
	W_{n, \max}^* \equiv \max_{1 \leq l \leq L} \left|{T}_{n,l}^{*{\scriptscriptstyle H_0}}\right|, 
\end{align}
with 
\[
{T}_{n,l}^{* {\scriptscriptstyle H_0}}  \equiv  \sqrt{n \over \tilde{s}_{h_n}^*(X_{0, l})}   
  {\langle \tilde{\beta}_{h_n}^*, X_{0,l} - \bar{X}^* \rangle},
\quad l=1,\ldots,L,  
\]
denoting the bootstrap rendition of the estimated projection quantities ${T}_{n,l}^{ {\scriptscriptstyle H_0}}$ from \eqref{eqn:Hostat} under $H_0$.  Above $\tilde{s}_{h_n}^*$ denotes estimated scaling, akin to $\hat{s}_{h_n}$, computed from the bootstrap sample $\{(X_i^*,\tilde{Y}_i^*)\}_{i=1}^n$.

The following result guarantees that, under the null hypothesis $H_0$, the distribution of test statistics $W_{n,L^2}$ and $W_{n, \max}$ in \eqref{eq_testing_stat} can be approximated by either bootstrap approach: enforcing $H_0$ as in \eqref{eq_testing_statTilde} or not as in  \eqref{eq_testing_statHat}.  

\begin{cor} \label{cor_test}
	Let $W_n$ denote a test statistic $W_{n,L^2}$ or $W_{n, \max}$ and let $W_n^*$ denote its paired bootstrap (PB) counterpart, computed   either as in \eqref{eq_testing_statHat} or  \eqref{eq_testing_statTilde}.  Under the assumptions of \autoref{thm_pb}, if the null hypothesis $H_0$ in \eqref{eq_null} holds, then  
	\begin{align*}
		& \sup_{w \in \R} |\pr^*(W_{n}^* \leq w|\xx_0) 
		- \pr ( W_{n} \leq w| \xx_0 )|  \xrightarrow{\pr} 0
		\quad \mbox{as $n\to \infty$}. 
	\end{align*} 
\end{cor}

While both implementations \eqref{eq_testing_statHat}-\eqref{eq_testing_statTilde} of PB are valid for testing, numerical studies suggest that enforcing the null hypothesis 
\eqref{eq_testing_statTilde} can have better performance in both  size and power. This is explored further in \autoref{ssec_4_2}.


\section{Simulation studies} \label{sec4}

\autoref{ssec_4_1} summarizes   numerical studies of the PB and other methods for calibrating confidence intervals for conditional response mean $\mu(X_0)$  in FLRMs.  A rule of thumb for selecting the tuning parameters $(k_n, h_n, g_n)$ in the bootstrap procedure is also examined. \autoref{ssec_4_2} then investigates the performance of the   bootstrap test   from \autoref{sec3} regarding projections.  {For simulation purposes, we consider the case when $\eo[X] = 0$ and $\eo[Y] = 0$ so that the conditional mean is equal to the projection of the slope function, i.e., $\mu(X_0) = \langle \beta, X_0 \rangle$. 
We hence focus on inferring the projection $\langle \beta, X_0 \rangle$ by using the original $\langle \hat{\beta}_{h_n}, X_0 \rangle$ and bootstrap $\langle \hat{\beta}_{h_n}^*, X_0 \rangle$ projection estimators with intercepts $\ap, \hat{\ap}_{h_n}, \hat{\ap}_{h_n}^*$  set to  zero in all following bootstrap or normal approximations for simplicity;
the sample means $\bar{X},\bar{X}^*$ used in the scaling terms are also adjusted to zero.
For PB, this is equivalent to using \autoref{cor_proj}  with  $\bar{X},\bar{X}^*$ being   zero.
}


\subsection{Performance of bootstrap intervals} \label{ssec_4_1}
 
To describe the data generation, we   independently simulate $n$ curves $\xx_n = \{X_i\}_{i=1}^n$ from a truncated Karhunen--Lo\`{e}ve expansion:
\begin{align} \label{eq_KL_trunc}
	X \ed \sum_{j=1}^J \sqrt{\g_j} \xi_j \phi_j
\end{align}
with $J=15$.  Above    $\{\phi_j:j =1, \dots, J \}$ denote the first $J$ of the Fourier basis functions $\{1, \sin(2\pi t), \cos(2\pi t), \dots\}$ on  $[0,1]$. 
The FPC scores are defined as $\xi_j = \xi W_j$, where $W_j \iid \nd(0,1)$ and $\xi$ follows  a {(normalized)} $t(\nu)$ distribution 
with chosen degrees  of freedom $\nu\in\{4,5,7,9,\infty\}$.
This entails that FPC scores are uncorrelated, but dependent. 
The eigengaps in \eqref{eq_KL_trunc} are defined with a polynomial decay rate involving a parameter $a>0$,
namely  $\gamma_{j}-\gamma_{j+1} =2j^{-a}$, $j \geq 1$
where $\g_1 = 2\sum_{j=1}^\infty j^{-a}$. 
Using the same basis functions,   the slope parameter is set to $\beta = \sum_{j=1}^{J} \beta_j \phi_j$, where  $\beta_j = 3j^{-b} W_{\beta,j}$   has a polynomial decay involving a rate parameter $b>0$ and the terms $W_{\beta, j}$ are fixed upon drawing these as iid from a distribution $\pr(W_{\beta, j}=1) = 1/2 = \pr(W_{\beta, j} = -1)$.  We consider various scenarios involving different polynomial rates and sample sizes: $a,b \in \{1.5, 2.5, 3.5, 4.5, 5.5\}$ and $n \in \{50, 200, 1000\}$.  For brevity, we report some representative numerical results here, though full results can be found in the   supplement \cite{supp}. All the function values are evaluated at 100 equally-spaced time grid points in $[0,1]$.  Response values $\{Y_i\}_{i=1}^n$ are  then generated through the FLRM \eqref{eq_model_obs} as follows.  
To consider both homoscedastic and heteroscedastic scenarios, 
errors $\e_i$ are generated to be either independent from or dependent on the regressors $X_i$.   For a given generated regressor  $X_i$, a dependent error $\e_i$ is simulated from a chi-square distribution $\chi^2(\nu(X_i)) - \nu(X_i)$ with   $\nu(X_i) \equiv\|X_i\|^2/2$ degrees of freedom. In this heteroscedastic case, the conditional variance of an error depends on the regressor value $X_i$, and the marginal variance of an error is $\var[\e_i] = \tr{\ga} = \sum_{j=1}^J \g_j$.
Due to the latter,  we also generate errors $\e_i$ with the same marginal variance, independently from regressors $X_i$, with a centered chi-square distribution $\chi^2(\nu) - \nu$ with $\nu \equiv \tr{\ga} / 2$ degrees of freedom in  homoscedastic cases.
The supplement \cite{supp} provides further results with other error distributions, which are qualitatively similar. In each simulation run, a regressor $X_0$ for mean $\mu(X_0)$ estimation is also generated by  \eqref{eq_KL_trunc}.

We consider both PB and naive PB implementations 
for computing two-sided 95\% intervals for $\mu(X_0) =  \langle \beta, X_0 \rangle$. 
In  the original data FPCR estimator $\hat{\beta}_{h_n}$ from \eqref{eq_fpcr} and estimated scaling $\hat{s}_n$ from \eqref{eq_sHat}, we varied the range of the 
truncation	  parameters  $h_n, g_n \in\{1, \dots, 15\}$ and  we set
$k_n = 2 [n^{1/v}]$ with $v = 2a+1+v_1$ for a small $v_1 = 0.1$
for consistent estimation of $\hat{\beta}_{k_n}$ in scaling \eqref{eq_sHat}  (cf.~Theorem~S3 of the supplement \cite{supp}).  
To recall, $k_n$ is used to reconstruct the residual as used in the scaling factor \eqref{eq_sHat}, $g_n$ is for constructing the bootstrap centering, and $h_n$ is the truncation used by the actual and bootstrap estimators (see \autoref{thm_pb}). 
For simplicity here, we focus on symmetrized intervals in the PB implementation involving bootstrap estimated studentization (e.g., $T_{n,\hat{s}^*}^*(X_0)$ in \eqref{eq_SnHatStar_std}) as well as a naive PB counterpart defined {in \eqref{eq_SnStarCheck}.}
{
	The symmetrized intervals are constructed by approximating the distribution of the absolute statistic $|T_n|$, instead of the original statistic $T_n$, with bootstrap (cf.~\cite{Hall88}).
	Namely, letting $\hat{c}$ be the $1-\ap/2$ quantile of the (bootstrap) distribution of $|T_n^*|$, where $T_n^*$ denotes the bootstrap statistic from either \eqref{eq_SnHatStar_std} or \eqref{eq_SnStarCheck},
	the symmetrized confidence interval for $\mu(X_0) =  \langle \beta, X_0 \rangle$ is constructed as $ \langle \hat{\beta}_n,X_0\rangle \pm \hat{c} \hat{s}_{h_n}(X_0) /\sqrt{n}$.
}
Further comparisons with other versions of PB (e.g., $T_{n,\hat{s}}^*(X_0)$ in \eqref{eq_SnHatStar})  or non-symmetrized intervals can be found in the supplement \cite{supp}, 
though estimated bootstrap studentization steps tend to induce the best performances.  
For comparison, we also consider intervals based on  normal approximations with estimated scaling $\hat{s}_n$ (\autoref{thm_clt}), residual bootstrap (RB)     (cf.~\cite{YDN23RB}), {and wild bootstrap (WB) \cite{HM91}.
The WB method shares similarities to the RB  up to the bootstrap error construction;
the latter involves a random weighting of residuals by the two-point distribution as described in \cite{HM91}.}
For each generated data set, bootstrap distributions are approximated by 1000 Monte Carlo resamples.

We also propose a rule of thumb for setting the tuning parameters based on simulations for all the parameter combinations.
We suggest to set $g_n=k_n$ and $h_n = [1.1k_n]$ being a slightly larger value than $g_n$; the value of $k_n$ can be selected in practice by cross-validation minimizing the prediction errors.
Our rule of thumb is found by considering all scenarios and truncation levels producing coverages of PB intervals within 1\% from the nominal level 95\%, and running linear regression of response $(h_n, g_n)$ on $k_n$.   
This rule targets to make appropriate choices of $(h_n,g_n)$, as   most critical to performance of PB, in relation to $k_n$.	
Setting $g_n=k_n$ aligns with the appropriate choices for the RB  \citep{YDN23RB}. 

\begin{figure}[b!]
	\centering
	\includegraphics[width=0.89\linewidth]{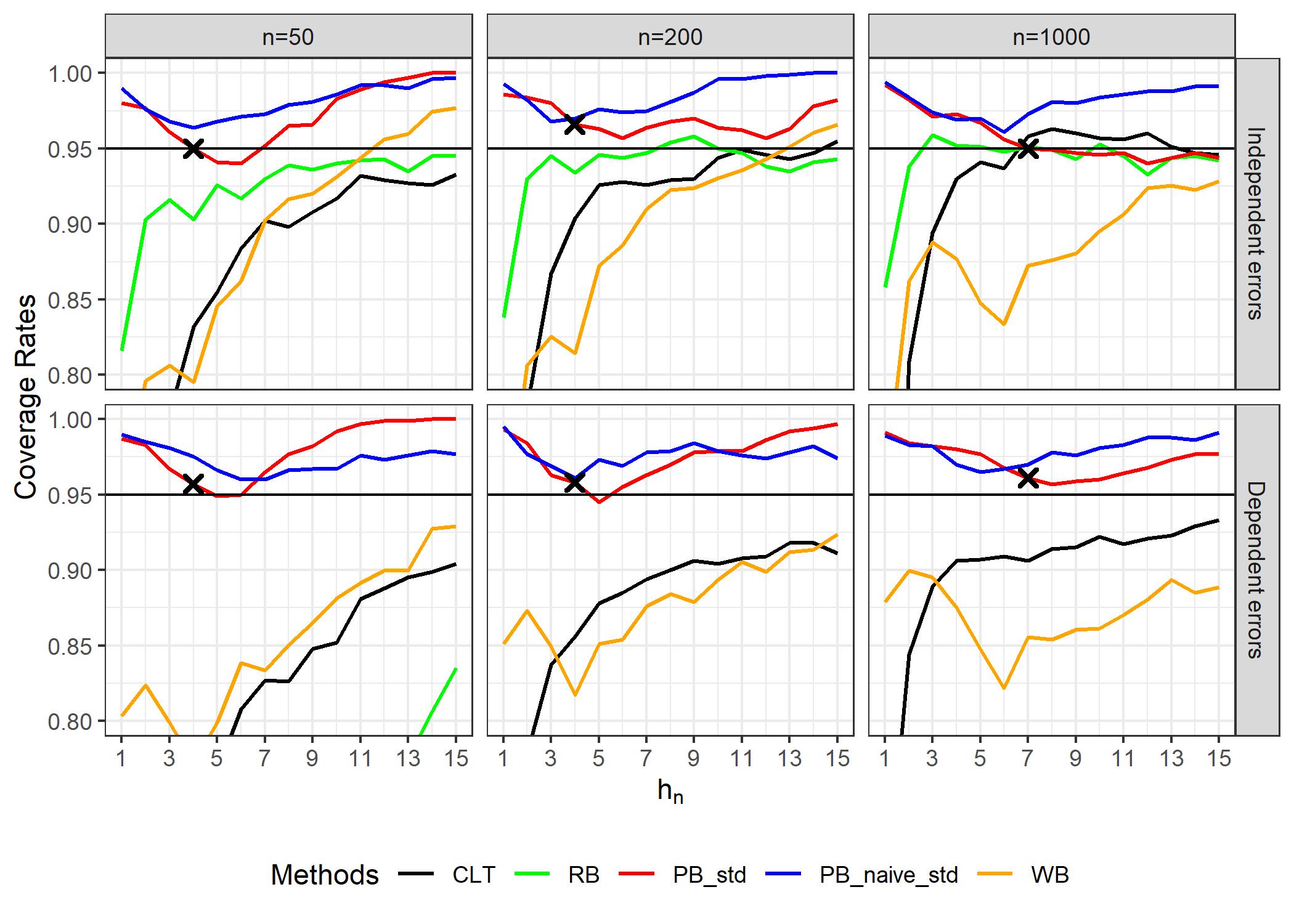}
	\caption{
		Empirical coverage rates of 95\% intervals for $  \mu(X_0)=\langle \beta, X_0 \rangle$ from CLT (black), RB (green), PB with   studentization (red), naive PB with studentization (blue),  and WB (tan)   over various truncations $h_n$ 
		when the decay rates for $\g_j - \g_{j+1}$ and $\beta_j$ are   $a=2.5$ and $b=5.5$ and the latent variable for the FPC scores is $\xi \sim t(5) / \sqrt{5/3}$. 
		For  errors dependent on regressors (lower panels), 
		the coverage curves of   CLT/RB intervals are cropped as these perform poorly.
		Crosses $\times$ indicate coverage rates with $h_n$ selected by the proposed rule.
	}
	\label{fig_cover1_crop}
\end{figure}

For each 95\% interval procedure for $\mu(X_0) = \langle \beta, X_0 \rangle$,  empirical coverages  were approximated by  1000 simulation runs over each data generating model and sample size.  
\autoref{fig_cover1_crop} displays observed coverage rates from different methods under a few selected scenarios when $a=2.5$, $b=5.5$ and $\xi \sim t(5) {/\sqrt{5/3}}$; see the supplement \cite{supp} for results over all scenarios.  
For clarity, the results in \autoref{fig_cover1_crop} focus on the case that $g_n=k_n$ for both PB and RB
while varying $h_n$.	
Coverages for the PB method under the proposed rule of thumb are indicated  using crosses in \autoref{fig_cover1_crop} for reference.    

As a first observation from \autoref{fig_cover1_crop},   the coverages from intervals based directly on normal approximation (CLT) exhibit sensitivity to the truncation level $h_n$ and also under coverage,  particularly when the sample size is small.  Under heteroscedasticity, both the CLT and RB methods perform quite poorly and lie at least partially outside of the charting regions in \autoref{fig_cover1_crop}. In fact, under this case of heteroscedasticity, RB is not asymptotically valid   and the coverages are quite low to the extent that coverage curves do not appear in the figure, even for large sample sizes $n=1000$.
{The WB, which is appropriate for heteroscedastic error structures, performs better than RB in this case, but still exhibits low coverage accuracy.}
In contrast, PB intervals perform much better under the heteroscedastic models. 
For independent errors, while RB assumes and uses the true model structure (homoscedasticity) and PB does not, the PB method has very similar performance to RB for large sample sizes ($n=1000$) and exhibits comparable performance for smaller samples ($n=50$ or $200$).  Our rule of thumb provides reasonable coverages in most cases for PB intervals.

\begin{figure}[b!]
	\centering
	\includegraphics[width=0.89\linewidth]{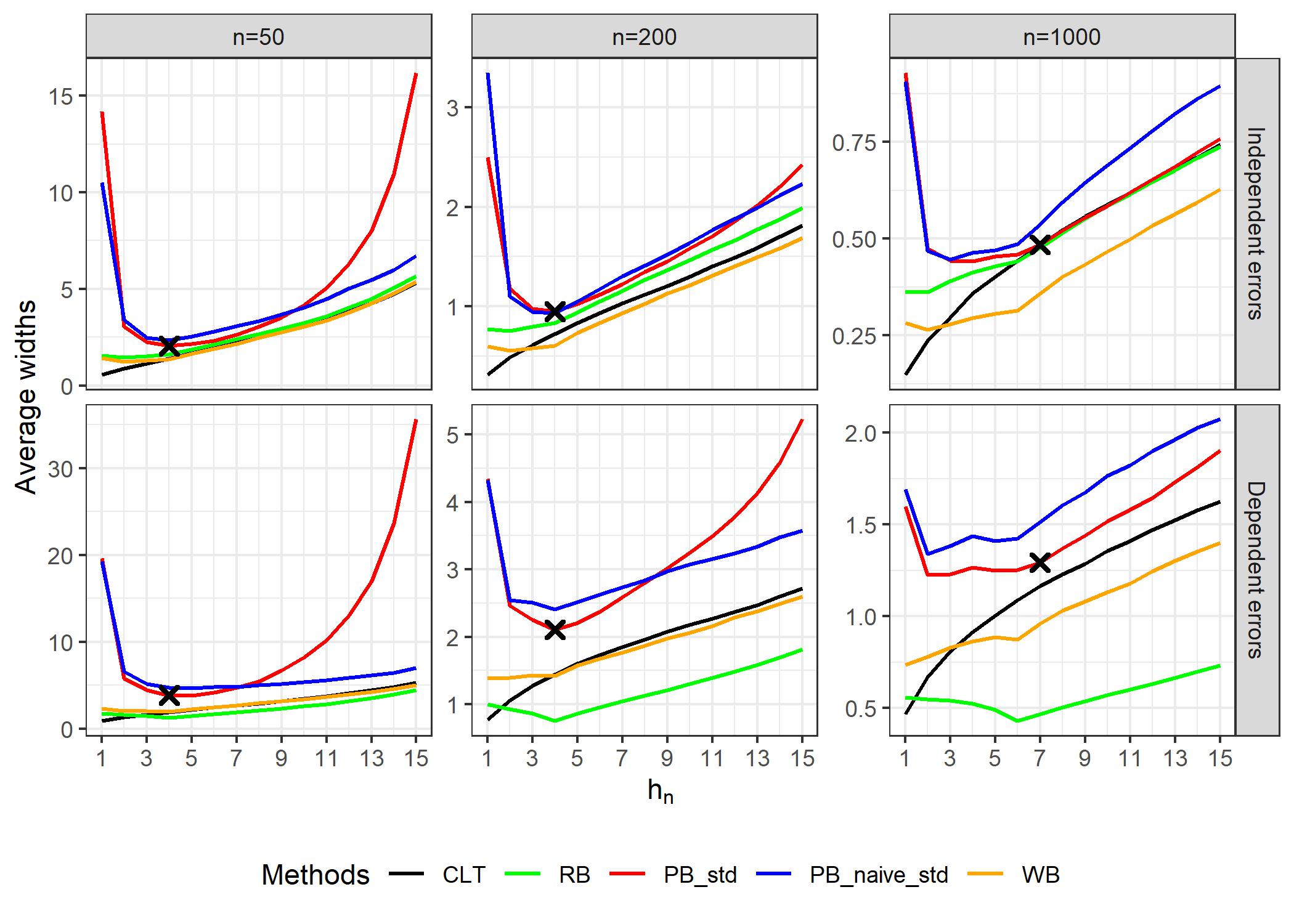}
	\caption{
		Average widths of  95\% intervals for $\mu(X_0)=\langle \beta, X_0 \rangle$ from   methods over different truncation $h_n$ 
		when the decay rates for $\g_j - \g_{j+1}$ and $\beta_j$ are  $a=2.5$ and $b=5.5$ and the latent variable for the FPC scores is $\xi \sim t(5) / \sqrt{5/3}$:  
		CLT (black), RB (green), WB(tan), PB with studentization (red), and naive PB with studentization (blue). Crosses $\times$ indicate coverage rates for PB with $h_n$ selected by the proposed rule.
	}
	\label{fig_width1}
\end{figure}

\autoref{fig_width1} displays the corresponding  average widths of intervals,  which generally increase with   $h_n$.
Importantly, this figure indicates that intervals from   RB  and CLT approximations are often overly narrow under heteroscedasticity, 
which relates to the low coverages in \autoref{fig_cover1_crop}.  Figures~\ref{fig_cover1_crop}-\ref{fig_width1} also demonstrate 
that our rule of thumb seems to suggest an optimal truncation $h_n$ 
in the sense that the corresponding intervals balance good coverage rates with lowest average widths.
Finally, while   the naive PB implementation is  asymptotically  invalid  in the sense of  \autoref{prop_pb_naive_fail},  the latter finding also suggests that  the bias in the naive PB should translate to over-coverage for symmetrized intervals in \autoref{fig_cover1_crop}.  
Even for large sample sizes $n=1000$,   naive PB intervals tend to over-cover   projections, while their average widths are larger than those from the proposed PB.  
Moreover, the coverages of naive PB intervals are unstable against the choice of truncation level $h_n$ 
while our modified PB produces stable coverages close to the nominal level for all $h_n \geq g_n$ and moderate to large sample sizes $n=200$ and 1000. 
The over-coverage problem in the naive PB also worsens 
as  truncation levels $h_n$ deviate from the case $h_n=g_n$.
This can be interpreted as   the construction bias from the naive bootstrap  negatively impacts this method, even as the sample size increases.

	To investigate the effect of  moments  for the regressor $X$ on interval performance, we   also varied the distribution of $\xi$ in \eqref{eq_KL_trunc} over different {normalized} $t(\nu)$   with  $\nu\in\{4,5,7,9,\infty\}$, where $t(4)$ is an example that does not satisfy \ref{condA2_4thmoment}.
	Figures in Section~S4.1  of the supplement \cite{supp} show that, 
	under  heteroscedasticity,  the proposed PB is fairly robust to   moments available  for  $X$.

\subsection{Performance of bootstrap tests of projections} \label{ssec_4_2}

Turning to the testing problem discussed in \autoref{sec3},
We next investigate the empirical rejection rates of the bootstrap tests 
of a null hypothesis \eqref{eq_null_proj} of constant conditional means, using bootstrap statistics 
 from \eqref{eq_testing_statHat} or \eqref{eq_testing_statTilde}.

The data generation for purposes of study   are generally the same as  considered in  \autoref{ssec_4_1} with $\xi \sim \nd(0,1)$ and $a=2.5$, with the exception that we modify the definition of the slope function $\beta$ to describe different hypotheses.  For testing,  the target predictors are considered as $\xx_0 \equiv \{ \phi_j \}_{j=1}^6$ based on the first six Fourier basis functions.
  Under the null hypothesis, the slope function is defined as $\beta^{H_0} \equiv \sum_{j > 6} W_{\beta, j} |\beta_j| \phi_j$,  and we wish to assess the orthogonality of $\beta$ to the subspace spanned by $\xx_0$ (i.e., $H_0:  \Pi_{\mathcal{X}_0}\beta=0$ or \eqref{eq_null_proj} with $\eo[X]=0$).  The true data-generating slope is defined as $\beta^{H_1} \equiv (1-p)\beta^{H_0} + p \sum_{j=1}^6 W_{\beta, j} |\beta_j| \phi_j$   in terms of a proportion $p\in\{0, 0.1, \dots, 0.9, 1\}$ for prescribing a sequence of alternative hypotheses;   
here  $|\beta_j| = c j^{-b}$ holds with $c = 50$ and $b = 3.5$, and the value $p=0$ renders the null hypothesis with the slope  $\beta^{H_0}$.

We consider bootstrap tests of $H_0:\Pi_{\mathcal{X}_0}\beta=0$ based on a nominal size 5\%.  
For each simulated dataset, 1000 bootstrap resamples are used to approximate the distribution of test statistics in \eqref{eq_testing_stat}.  Truncation parameters $h_n$ and $g_n$ are again selected by the rule of thumb suggested in \autoref{ssec_4_1} based on  $k_n$.    
Using 1000 simulation runs for each data generation scenario (level of $p$) and sample size $n$,  we compute rejection rates by the proportion of times that an original test statistic exceeds the 95th percentile of bootstrap test distribution.   
The supplement \cite{supp} contains more details and findings over 
different sample sizes $n\in\{50, 200, 400, 600, 800, 1000\}$ 
as well as both test statistic forms from \eqref{eq_testing_stat}; we present results for $n=50$ here with  
maximum or $L_\infty$ statistic form $W_{n,\max}$, as other results are qualitatively similar.

\begin{figure}[b!]
	\centering
	\includegraphics[width=0.59\linewidth]{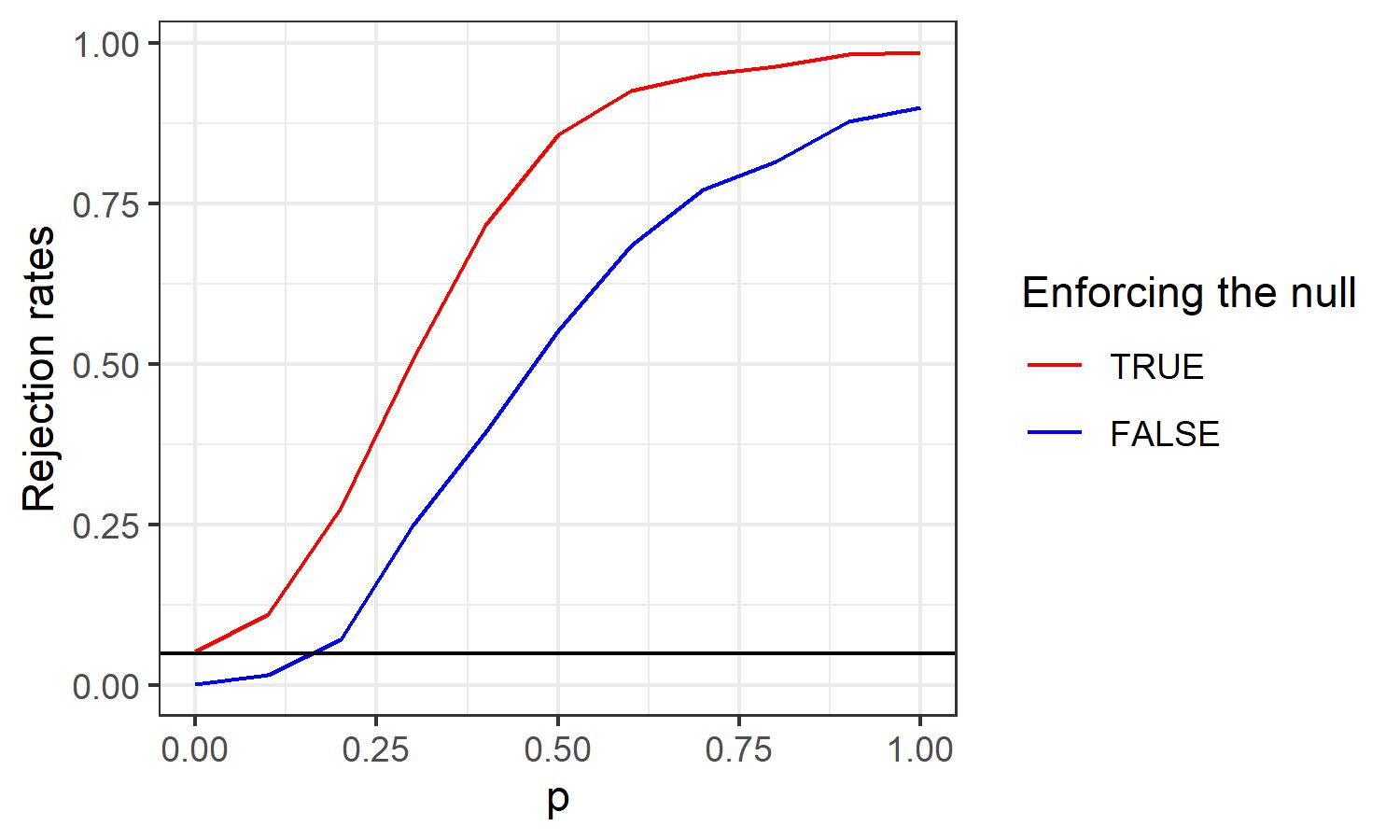}
	\caption{Empirical rejection rates (when $n=50$) of the bootstrap testing procedure as the degree/proportion $p\in\{0, 0.25, 0.50, 0.75, 1\}$ of the alternative increases (only $p=0$ corresponds to a true null hypothesis). The test may enforce the null hypothesis (red) or not (blue) in the bootstrap. The black horizontal line represents the nominal size 5\%.}
	\label{fig_rejRate}
\end{figure}

The resulting empirical rejection rates are summarized in \autoref{fig_rejRate}. 
As perhaps expected, the power of the test  increases with the degree $p$ of how much the null hypothesis is violated, whether enforcing the null hypothesis in bootstrap by \eqref{eq_testing_statTilde},   or not   by \eqref{eq_testing_statHat}.  However, enforcing the null hypothesis maintains size better (i.e., when $p=0$), which then also leads to slightly better power here.   Another advantage to bootstrap enforcement of the null hypothesis
is less sensitivity to choices of truncation parameters $h_n,g_n$. Results in the supplement \cite{supp} indicate that honoring the null hypothesis in bootstrap  typically ensures good performance in testing as  truncations $h_n,g_n$ are varied, which  is not equally true for the bootstrap version that does not enforce the null hypothesis.

\section{Real data analysis} \label{sec5}

\begin{figure}[b!]
	\centering
	\includegraphics[width=\linewidth]{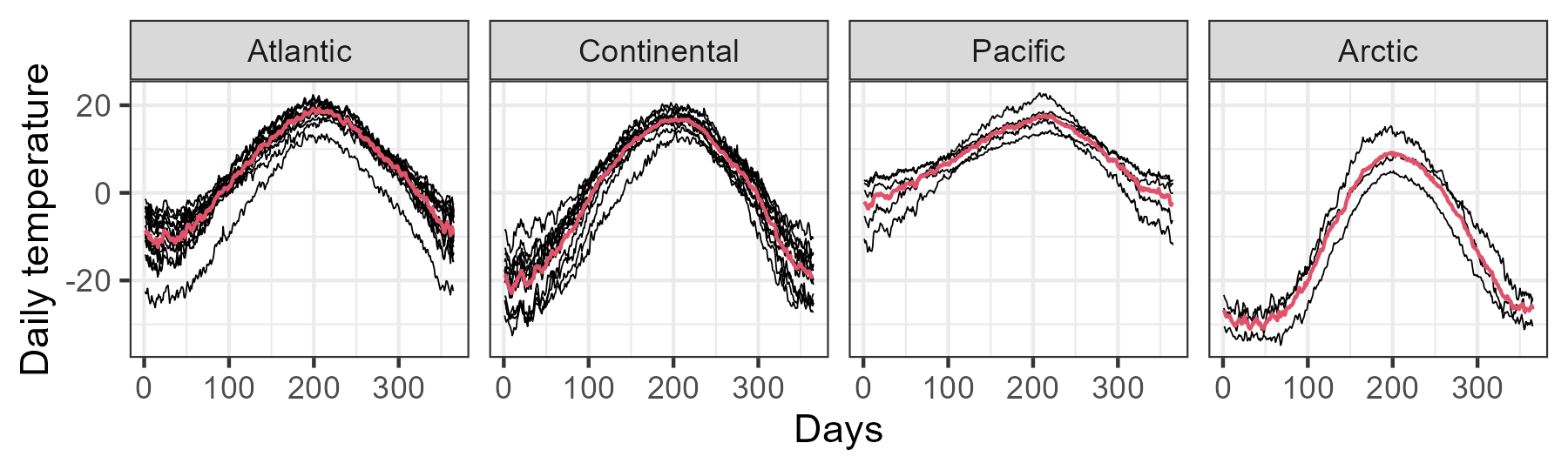}
	\caption{Daily temperature curves of locations in four different regions. Each black curve corresponds to the averaged in one location over 1964 to 1990, and the regional average curves are denoted in bold pink lines.}
	\label{fig_rda_cw_data}
\end{figure}

Bootstrap   intervals and tests   are demonstrated in application to   Canadian weather data.
We analyze the Canadian weather dataset from the R package \texttt{fda} consisting of daily temperature and precipitation at 35 different locations in Canada \cite[cf.][]{RS05}.
The regressor $X_i$ is the daily temperature on each day averaged over 1960 to 1994, and the response $Y_i$ is the log of total annual precipitation with base 10. 
Pairs  of temperature  and precipitation  are recorded at $n=35$ weather stations. 
The regressor curves $\{ X_i \}_{i=1}^n$ are displayed in \autoref{fig_rda_cw_data}, where the thicker lines represent the average for four different regions in Canada, namely, Atlantic, Continental, Pacific, and Arctic regions.
The new predictors $\xx_0 \equiv \{X_{0,l}\}_{l=1}^4$ under consideration for bootstrap inference are selected as these four average curves in each region as  illustration.
{
We then focus  inference of centered projections $\langle \beta, X_{0,l}-\eo[X] \rangle$ (cf.~\autoref{cor_proj}),
which essentially represents the conditional mean $\mu_0(X_{0,l}) =\eo[Y] + \langle \beta, X_{0,l}-\eo[X] \rangle$
at each regressor $X_{0,l}$ up to a common effect $\eo[Y]$   not involving  the slope $\beta$.  }

Each weather station is located in one of the four regions, where each region exhibits a different pattern as shown in \autoref{fig_rda_cw_data}. 
This leads us to suspect the existence of different conditional variance of errors in FLRM \eqref{eq_model_intercept}.
To investigate the heteroscedasticity, we   estimate the variance from residuals for each region as $\hat{\s}_{r, k_n}^2 =  n_r^{-1} \sum_{i \in \ii_r} (Y_i -\bar{Y} - \langle \hat{\beta}_{k_n}, X_i-\bar{X} \rangle )^2 $, where   $\ii_r$ and  $n_r$, respectively, denote the index set of and the number of location in the $r$th region. 
Here, the estimator $\hat{\beta}_{k_n}$ used for computing residuals is constructed from the combined data $ \{(X_i, Y_i)\}_{i=1}^n$ over all four regions.
As shown in \autoref{fig_rda_cw_estSD}, homoscedastic error models seems implausible for this dataset. 
A similar conclusion is deduced from the residual plots given in Section~S5.1
of the supplement.

{By applying different bootstrap methods,  symmetrized  $95\%$   confidence intervals
for each projection $\langle \beta, X_{0, l} - \eo[X] \rangle$ are given in \autoref{tb_cw}. Based on PB, such intervals 
are   $\langle \hat{\beta}_n, X_{0,l}-\bar{X} \rangle \pm \hat{c} \hat{s}_{h_n}(X_0) /\sqrt{n}$, with $\hat{c}$ denoting the 95th percentile 
of the absolute bootstrap statistic 
$|T_{n,proj}^*|$ from  \eqref{eq_stat_proj_bts} (denoted as PB\_std) or from
\eqref{eq_stat_proj_bts1}
(denoted as PB);  recall PB\_std uses estimated bootstrap studentization, while PB does not, for improved accuracy.} Here, the less consequential tuning parameter $k_n=2$ was selected via repeated cross-validation, which minimizes prediction errors over estimates from $\hat{\beta}_{k_n}$, while $h_n = 2$ and $g_n = 2$ were then chosen by the rule of thumb suggested in \autoref{ssec_4_1}. 
The supplement \cite{supp} provides further results with different tuning parameter choices.
As expected under possible heteroscedasticity and shown in \autoref{tb_cw}, the results for residual bootstrap (RB) are quite different from those for PB.
This distinction  is also seen from a comparison of interval lengths in \autoref{tb_cw}.
Compared to an overall average $\eo[X]$, the Pacific region has the highest range of annual precipitation, while the Continental region exhibits less precipitation, with relatively narrow widths for both regions.
The annual precipitations of the Atlantic and Arctic regions are, respectively, either higher or lower than the overall average, but with wider widths. 

We may also apply our bootstrap testing procedure  for assessing the null hypothesis $H_0$ from \eqref{eq_null} of the equality of means $\mu(X_{0,l})$ across the four regions.  
The corresponding p-values are given in \autoref{tb_cw_testing} for 
test statistics \eqref{eq_testing_statTilde} based on PB approximations that enforce $H_0$.  The tests    strongly suggest that the true  rainfall mean responses $\{ \mu(X_{0,l}) \}_{l=1}^4$ are not equal across regions and 
cannot be simultaneously equal to a common mean response $\mu(\eo[X])$ at the global mean curve.  
  This finding supports  the region-wise intervals in  \autoref{tb_cw}.

\begin{table}[t!]
	
	\centering
	\caption{
		$95\%$ symmetrized confidence intervals for centered projections $\{\langle \beta, X_{0, l} - \eo[X] \rangle\}_{l=1}^4$ from RB, PB, and PB\_std  for Canadian weather dataset.
		The ratios of widths of RB intervals to widths of either PB or PB\_std intervals are given in parentheses.
	} 
	\begin{tabular}{c|c|c|c}
		\hline
		& RB & PB & PB\_std
		\\ \hline 
		Atlantic 		& [~0.06, ~0.11] & [~0.00, ~0.17] (3.20) & [-0.05, ~0.22] (5.11)
		\\ Continental	& [-0.19, -0.08] & [-0.22, -0.05] (1.53) & [-0.25, -0.03] (1.96)
		\\ Pacific		& [~0.18, ~0.36] & [~0.17, ~0.38] (1.18) & [~0.14, ~0.41] (1.50)
		\\ Artic		& [-0.49, -0.18] & [-0.57, -0.10] (1.53) & [-0.57, -0.10] (1.54)
		\\ \hline
	\end{tabular}
	\label{tb_cw}
\end{table}

\begin{table}[b!]
	\centering
	\caption{P-values for   bootstrap tests to assess null hypothesis of constant rainfall means $\mu(X_{0,l})= \mu(\eo[X])$ across four regions    for Canadian weather dataset.}
	\begin{tabular}{c|c}
		L2 Statistic & Max Statistic \\  \hline
		  0.004   & 0.012 \\ 
		\hline
	\end{tabular}
	\label{tb_cw_testing}
\end{table}


\section{Concluding remarks} \label{sec6}

We have developed a paired bootstrap (PB)  for inference in functional linear regression models  (FLRMs) with general heteroscedastic errors.  
As a preliminary result, a central limit theorem  under heteroscedasticity was established for estimated conditional mean $\mu(X_0) \equiv \ap + \langle \beta, X_0 \rangle$ given new predictor $X_0$, 
along with appropriate scaling $s_{h_n}(X_0)$ for self-normalization.  
Further, the estimated conditional mean $\hat{\mu}_{h_n}(X_0) \equiv \hat{\ap}_{h_n} + \langle \hat{\beta}_{h_n}, X_0 \rangle$ given new predictor $X_0$ based on  a functional principal component regression estimator $\hat{\beta}_{h_n}$  can be successfully approximated by the PB for improved inference in finite samples.  
As such estimators  $\hat{\beta}_{h_n}$ in FLRMs involve truncation parameters $h_n$, a modified PB   estimator was  proposed  to allow valid distributional approximations with the  greatest possible flexibility in such truncation parameters for bootstrap. In contrast,   a naive implementation of PB (i.e., adapted directly from standard multiple regression, e.g., \cite{free81}) can be shown to have less validity in application and becomes viable
only for a much narrower configuration of truncation parameters.    
The PB approach was also adapted to formulate new tests for assessing the equality of conditional means at pre-selected regressor curves $\{X_{0,l}\}_{l=1}^L$, which can boil down to the orthogonality  of the slope function $\beta$ to the subspace spanned by $\{X_{0,l}\}_{l=1}^L$ under zero-mean assumption. 
Numerical studies showed that the existing residual bootstrap can fail  under heterocedasticity,  while PB can perform well in this context for interval estimation as well as for testing.
We suggested also a PB based on estimated bootstrap studentization steps and provided a rule of thumb for selecting the two main tuning parameters (truncation levels) involved in the PB. 

Other  directions are possible for bootstrap inference in the FLRMs, which require further study. 
A wild bootstrap method \citep{DBZ17, HM91, mammen93} might also be formally developed for inference in FLRMs with heteroscedastic errors. A PB approach can have more accuracy than  wild bootstrap in some settings (cf.~\autoref{ssec_4_1}), 
but wild bootstrap may be computationally faster and more scalable for large data in its resampling scheme;
further investigation is needed for FLRMs in particular.   
It can also be of  interest to  extend resampling inference in other functional linear models such as FLRMs with functional response \citep{CM13} or generalized functional linear models \citep{MS05}.  {Finally, as described in Remark~6, further investigation is required to extend bootstrap inference  to cases where functional regressors may be irregularly or sparsely observed (cf.~\cite{ZYZ23}).}


\begin{appendix}

\section{On establishing the CLT}  \label{app1}
We  briefly outline of the proof of the CLT in \autoref{thm_clt}; more technical details are provided 
in the supplement \cite{supp}.

\begin{proof}[Proof of \autoref{thm_clt}]
	{
	It suffices to show the weak convergence of the projection estimator $\langle \hat{\beta}_{h_n}, X_0 \rangle $ under zero-mean assumptions. 
	To explain, the conditional mean difference is decomposed as 
	\begin{align} \label{eq_thm_clt_condMean}
		\hat{\mu}_{h_n}(X_0) - \mu(X_0)
		= \bar{Y} - \eo[Y] - \langle \hat{\beta}_{h_n}, \bar{X} - \eo[X] \rangle  + \langle \hat{\beta}_{h_n} - \beta, X_0 - \eo[X] \rangle .
	\end{align}
	Using that $\sqrt{n}(\bar{Y} - \eo[Y]) = O_\pr(1)$ by classical CLT \cite[cf.][]{bill95}, that $\sqrt{n}(\bar{X} - \eo[X]) = O_\pr(1)$ by   CLT in general Hilbert space \cite[cf.][]{HE15}, and that $\|\hat{\beta}_{h_n} - \beta\| = o_\pr(1)$ by the consistency of the FPCR estimator in the supplement \cite{supp},
	the first two terms in \eqref{eq_thm_clt_condMean} are asymptotically ignorable with scaling, i.e., 
	\begin{align*}
		\sqrt{n \over s_{h_n}(X_0)} (\bar{Y} - \eo[Y] - \langle \hat{\beta}_{h_n}, \bar{X} - \eo[X] \rangle) = O_\pr(1)s_{h_n}(X_0)^{-1} = O_\pr(h_n^{-1}).
	\end{align*}
	Because centering by the population means $\eo[Y]$ and $\eo[X]$ does not affect the construction of the estimator $\hat{\beta}_{h_n}$,
	it is enough to show the CLT for the projection with centered data $\{(\check{X}_i, \check{Y}_i)\}_{i=1}^n$ where $\check{X}_i \equiv X_i - \eo[X]$ and $\check{Y}_i \equiv Y_i - \eo[Y]$.
	Hence, without loss of generality, we may assume that $\eo[X] = 0$ and $\eo[Y] = 0$ so that $\ap = 0$, 
	and establish the CLT for   $\langle \hat{\beta}_{h_n} - \beta, X_0   \rangle$.
	}

	The proof
	uses the following bias-variance decomposition of the functional principal component estimator $\hat{\beta}_{h_n}$ with respect to the true slope parameter $\beta$:\begin{align} \label{eq_decomp}
		\hat{\beta}_{h_n} - \beta
		= b_n + \ga_{h_n}^{-1} U_n,
	\end{align}
	where, upon scaling  $\sqrt{n / s_{h_n}(X_0)}$, the quantity $\ga_{h_n}^{-1} U_n$ determines the normal limit while a remainder/bias term $b_n \equiv \hat{\beta}_{h_n} - \beta - \ga_{h_n}^{-1} U_n$ converges to zero in probability.
	Above  $U_n \equiv n^{-1} \sum_{i=1}^n (X_i - \bar{X}) (\e_i - \bar{\e})$ represents the cross-covariance between the regressors $\xx_n \equiv \{X_i\}_{i=1}^n$ and the errors $\{\e_i\}_{i=1}^n$, with $\bar{X} \equiv n^{-1} \sum_{i=1}^n X_i$ and $\bar{\e} \equiv n^{-1} \sum_{i=1}^n \e_i$, and further $\ga_{h_n}^{-1} \equiv \sum_{j=1}^{h_n} \g_j^{-1} \pi_j$ denotes a truncated version of the inverse covariance operator $\ga^{-1} \equiv \sum_{j=1}^\infty \g_j^{-1} \pi_j$ with $\pi_j \equiv \phi_j \otimes \phi_j$ for integer $j \geq 1$.
	The supplement \cite{supp} shows  that, as $n\to\infty$,
	\begin{equation}
		\label{eqn:cltbias}
		\pr \left( \sqrt{n \over s_{h_n}(X_0)} | \langle b_n, X_0 \rangle | > \eta \Big| X_0 \right) \xrightarrow{\pr} 0
	\end{equation}
	holds for each $\eta>0$.     The distributional convergence of the   term $\ga_{h_n}^{-1} U_n$ is stated in the following proposition,
	where the proof is deferred to the supplement \cite{supp}.

	\begin{prop} \label{prop_clt_var_hetero}
		Suppose  that Conditions \ref{condA5_eg_rate}-\ref{condA7_4thmoment_xe} hold
  and $n^{-\dt/2} h_n^{(2+\dt)/2} \to 0$ for $\dt \in (0,2]$ in Condition~\ref{condA7_4thmoment_xe}, then 
  as $n\to\infty$ 
		\begin{align*}
			\sup_{y \in \R} \left| \pr \left( \sqrt{n \over s_{h_n}(X_0)} \langle \ga_{h_n}^{-1} U_n, X_0 \rangle \leq y \Big| X_0 \right) - \Phi(y) \right| \xrightarrow{\pr} 0.
		\end{align*}	 
	\end{prop}
	
	\autoref{thm_clt} then follows from (\ref{eqn:cltbias}) and 
	\autoref{prop_clt_var_hetero} under the decomposition \eqref{eq_decomp};  see also  Propositions~S1-S3 in the supplement \cite{supp}.  
\end{proof}

\section{On proofs for the paired bootstrap}  \label{app2}
We sketch the proofs of   \autoref{thm_pb}   and \autoref{prop_pb_naive_fail}; further details appear
in the supplement \cite{supp}.

\begin{proof}[Proof of \autoref{thm_pb}]
	
	{
	Similarly to \autoref{thm_clt}, 
	it suffices to consider the bootstrap   projection $\langle \hat{\beta}_{h_n}^* - \hat{\beta}_{g_n}, X_0 - \eo[X] \rangle$ with $\eo[X]=0$, as the difference from the bootstrap conditional mean 
	$  \hat{\mu}_{h_n}^*(X_0) - \hat{\mu}_{g_n}(X_0)
		- \langle \hat{\beta}_{h_n}^* - \hat{\beta}_{g_n}, X_0 - \eo[X] \rangle = 
		  \bar{Y}^* - \bar{Y}  - \langle \hat{\beta}_{h_n}^*, \bar{X}^* - \bar{X} \rangle  - \langle \hat{\beta}_{h_n}^* - \hat{\beta}_{g_n}, \bar{X} - \eo[X] \rangle$ is analogously negligible.  
	}
	
	To show  bootstrap consistency, we consider a bootstrap-level decomposition, similar to \eqref{eq_decomp}, as
	\begin{equation}
		\label{eqn:dn}
		\hat{\beta}_{h_n}^* - \hat{\beta}_{g_n} = b_n^* + \ga_{h_n}^{-1} (U_n^* - \hat{U}_{n,g_n})
	\end{equation}
	where $b_n^*$ is a bias term, $\hat{U}_{n,g_n}$ is the bias correction  from \eqref{eq_UnHatgn}, and $U_n^* \equiv n^{-1} \sum_{i=1}^n (X_i^* - \bar{X}^*) (\e_{i,g_n}^* - \bar{\e}^*_{g_n})$ denotes the sample cross covariance between the bootstrap regressors $\{X_i^*\}_{i=1}^n$ and the bootrstrap errors $\{\e_{i,g_n}^*\}_{i=1}^n$, where $\bar{X}^* \equiv n^{-1} \sum_{i=1}^n X_i^*$ and $\bar{\e}^*_{g_n} \equiv n^{-1} \sum_{i=1}^n \e_{i,g_n}^*$ from $\e_{i,g_n}^* \equiv Y_i^* - \langle \hat{\beta}_{g_n}, X_i^* \rangle$. \autoref{prop_pb_var} shows that, upon scaling, the distribution of $\ga_{h_n}^{-1} (U_n^* - \hat{U}_{n,g_n})$ under   bootstrap probability $\pr^*(\cdot|X_0)$ converges to a standard normal distribution.

	\begin{prop} \label{prop_pb_var}
		Suppose that Conditions \ref{condA1_model}-\ref{condA7_4thmoment_xe}, $\|\hat{\beta}_{g_n} - \beta\| \xrightarrow{\pr} 0$, and for $\dt \in (0,2]$ in Condition~\ref{condA7_4thmoment_xe}, we have $\eo[\|X\|^{4+2\dt}] <\infty$ and as $n\to\infty$, $n^{-\dt/2} h_n^{\dt/2} \sum_{j=1}^{h_n} \ld_j^{-(2+\dt)/2} = O(1)$ (which implies $n^{-\dt/2}h_n^{(2+\dt)/2} = o(1)$).
		Then, as $n\to\infty$, 
		\begin{align*}
			\sup_{y \in \R} \left| \pr^* \left( \sqrt{n \over s_{h_n}(X_0)} \langle \ga_{h_n}^{-1} (U_n^* - \hat{U}_{n, g_n}), X_0 \rangle \leq y \Big| X_0 \right) - \Phi(y) \right| \xrightarrow{\pr} 0.
		\end{align*}
	\end{prop}
	The supplement \cite{supp} then establishes that a scaled projection involving    
	$b_n^* \equiv \hat{\beta}_{h_n}^* - \hat{\beta}_{g_n} - \ga_{h_n}^{-1} (U_n^* - \hat{U}_{n,g_n})$ from (\ref{eqn:dn}) converges to zero in  bootstrap probability $\pr^*(\cdot|X_0)$, namely, 
	\begin{equation}\label{eqn:dn4}
		\pr^* \left( \sqrt{n \over s_{h_n}(X_0)} |\langle b_n^*, X_0 \rangle| > \eta \Big| X_0 \right) \xrightarrow{\pr} 0,
	\end{equation}
	as $n\to\infty$, for each $\eta>0$.  Using a subsequence argument (cf.~ \cite{bill95},~Theorem~20.5) for bootstrap distributions along with Slutsky's theorem,  \autoref{thm_pb} then follows from (\ref{eqn:dn4}) in combination with  \autoref{prop_pb_var} and (\ref{eqn:dn}); see also Propositions~S5-S9 in the supplement \cite{supp}. 
	
\end{proof}

\begin{proof}[Proof of \autoref{prop_pb_var}]
	We write $Z_{i,n}^* = \langle X_i^*\e_{i,g_n}^* - \tilde{U}_{n, g_n}, \ga_{h_n}^{-1} X_0 \rangle$ with bootstrap errors $\e_{i,g_n}^* \equiv Y_i^* - \langle \hat{\beta}_{g_n}, X_i^* \rangle$ and $\tilde{U}_{n, g_n} \equiv n^{-1} \sum_{i=1}^n X_i \hat{\e}_{i,g_n}$ so that $\eo^*[Z_{i,n}^*|X_0] = 0$ and 
	\begin{align}
		& \sqrt{n \over s_{h_n}(X_0)} \langle \ga_{h_n}^{-1} (U_n^* - \hat{U}_{n, g_n}), X_0 \rangle \nonumber
		\\=& \{n s_{h_n}(X_0)\}^{-1/2} \sum_{i=1}^n Z_{i,n}^*
		- \sqrt{n \over s_{h_n}(X_0)} \langle \bar{X}^*\bar{\e}_{g_n}^* - \bar{X}\bar{\hat{\e}}_{g_n}, \ga_{h_n}^{-1} X_0 \rangle, \label{PBeqBTSvar1}
	\end{align}
	where $\bar{\hat{\e}}_{g_n} \equiv n^{-1} \sum_{i=1}^n \hat{\e}_{i, g_n}$ with $\hat{\e}_{i, g_n} \equiv Y_i - \langle \hat{\beta}_{g_n}, X_i \rangle$.
        The second term in \eqref{PBeqBTSvar1}  negligible due to the following convergence:
	\begin{align} \label{PBlemBTSvar2_eq}
		\eo^* \left[ \left( \sqrt{n \over s_{h_n}(X_0)} \langle \bar{X}^*(\overline{\e^*})_{n, g_n} - \bar{X}(\bar{\hat{\e}})_{g_n}, \ga_{h_n}^{-1} X_0 \rangle \right)^2  \Big| X_0 \right] =o_\pr(1),
	\end{align}
	which follows using that, from Lemma~S35 in the supplement \cite{supp},
	\begin{align*}
		& \eo^* \left[ {n \over s_{h_n}(X_0)} \langle \bar{X}^*(\overline{\e^*})_{n, g_n} - \bar{X}(\bar{\hat{\e}})_{g_n}, \ga_{h_n}^{-1} X_0 \rangle^2  \Big| X_0 \right]
		\\ \leq & \{h_n s_{h_n}(X_0)^{-1}\} (n\eo^*[\|\bar{X}^*(\overline{\e^*})_{n, g_n} - \bar{X}(\bar{\hat{\e}})_{g_n}\|^2]) (h_n^{-1}\|\ga_{h_n}^{-1} X_0\|^2)
		\\ = &  \;O_\pr \left( n^{-1} h_n^{-1} \sum_{j=1}^{h_n} \g_j^{-1} \right),
	\end{align*}
	where the last big $O_\pr$ term converges to zero under Condition \ref{condA5_eg_rate}.
 
	To deal with the first term in \eqref{PBeqBTSvar1}, conditional on $X_0$, define the bootstrap variance $\hat{v}_n^2 \equiv \sum_{i=1}^n \eo^*[Z_{i,n}^{*2}|X_0]$ and a bootstrap version of the Lyapounov condition as $\hat{\kll}_n \equiv \hat{v}_n^{-(2+\dt)} \sum_{i=1}^n \eo^*[(Z_{i,n}^*)^{2+\dt}|X_0]$ for $\dt \in (0, 2]$.
	To verify a Lyapounov condition 
	\begin{align} \label{PBeqBTSvarLyapounov}
		\hat{\kll}_n \xrightarrow{\pr} 0,
	\end{align}
  we will use the following observation:
	\begin{align}
		\left| {n^{-1}\hat{v}_n^2 \over s_{h_n}(X_0)} - 1 \right| \xrightarrow{\pr} 0. \label{PBlemBTSvar1_eq}
	\end{align}
	To see the convergence in \eqref{PBlemBTSvar1_eq}, note that 
	\begin{align*}
		\eo^*[Z_{i,n}^{*2}|X_0]
		& = \eo^* [\langle X_i^*\e_{i, g_n}^* - \eo^*[X_i^*\e_{i, g_n}^*], \ga_{h_n}^{-1} X_0 \rangle^2|X_0]
		\\& = \langle \eo^* [(X_i^*\e_{i, g_n}^* - \eo^*[X_i^*\e_{i, g_n}^*])^{\otimes 2}] \ga_{h_n}^{-1} X_0, \ga_{h_n}^{-1}X_0 \rangle
	\end{align*}
	with $\eo^* [(X_i^*\e_{i, g_n}^* - \eo^*[X_i^*\e_{i, g_n}^*])^{\otimes 2}] = \eo^* [(X_i^*\e_{i, g_n}^*)^{\otimes 2}]  - (\eo^*[X_i^*\e_{i, g_n}^*])^{\otimes 2}$.
	We then find that  
	\begin{align*}
		\eo^* [(X_i^*\e_{i, g_n}^* - \eo^*[X_i^*\e_{i, g_n}^*])^{\otimes 2}]
		= \hat{\Ld}_{n, g_n},
	\end{align*}
	which implies that
	\begin{align*}
		\eo^*[Z_{i,n}^{*2}|X_0] 
		&= \langle \hat{\Ld}_{n, g_n} \ga_{h_n}^{-1} X_0, \ga_{h_n}^{-1}X_0 \rangle
		  = \langle (\hat{\Ld}_{n, g_n} - \Ld) \ga_{h_n}^{-1} X_0, \ga_{h_n}^{-1}X_0 \rangle + s_{h_n}(X_0),
	\end{align*}
	and hence, 
	\begin{align*}
		\left| {n^{-1}\hat{v}_n^2 \over s_{h_n}(X_0)} - 1 \right|
		= s_{h_n}(X_0)^{-1} |\langle (\hat{\Ld}_{n, g_n} - \Ld) \ga_{h_n}^{-1} X_0, \ga_{h_n}^{-1}X_0 \rangle|.
	\end{align*}
	The convergence in \eqref{PBlemBTSvar1_eq} now follows from Lemma~S30 in the supplement \cite{supp}.

	To prove \eqref{PBeqBTSvarLyapounov}, the Lyapunov term $\kll_n$ is expanded as
	\begin{align}
		\hat{\kll}_n 
		& = \hat{v}_n^{-(2+\dt)} \sum_{i=1}^n \eo^*[|Z_{i,n}^*|^{2+\dt}|X_0] \nonumber
		\\& \leq \hat{v}_n^{-(2+\dt)} \sum_{i=1}^n \eo^*[\|\Ld_{h_n}^{-1/2}X_i^* \e_{i, g_n}^*\|^{2+\dt} | X_0] \|\Ld_{h_n}^{1/2} \ga_{h_n}^{-1} X_0\|^{2+\dt} \label{PBeqBTSvar2}
		\\&= \left( \|\Ld_{h_n}^{1/2} \ga_{h_n}^{-1} X_0\|^2 \over n^{-1}\hat{v}_n^2 \right)^{(2+\dt)/2} n^{-\dt/2} \eo^* \left[ n^{-1} \sum_{i=1}^n \|\Ld_{h_n}^{-1/2}X_i^* \e_{i, g_n}^*\|^{2+\dt} \Big| X_0 \right]; \nonumber
	\end{align}
	here, we have
	\begin{align*}
		{\|\Ld_{h_n}^{-1/2} \ga_{h_n}^{-1}X_0\|^2 \over n^{-1} \hat{v}_n^2}
		\leq {s_{h_n}(X_0) \over n^{-1} \hat{v}_n^2} 
		= 1 + o_\pr(1) 
		= O_\pr(1)
	\end{align*}
	since $\|\Ld_{h_n}^{-1/2} \ga_{h_n}^{-1}X_0\|^2 = \langle \Ld_{h_n} \ga_{h_n}^{-1} X_0, \ga_{h_n}^{-1} X_0 \rangle \leq \langle \Ld \ga_{h_n}^{-1} X_0, \ga_{h_n}^{-1} X_0 \rangle = s_{h_n}(X_0)$
	and the last upper bound $O_\pr(1)$ is obtained due to \eqref{PBlemBTSvar1_eq}.
	The latter term in \eqref{PBeqBTSvar2} is bounded as
\begin{align}
		& n^{-\dt/2} \eo^* \left[ n^{-1} \sum_{i=1}^n \|\Ld_{h_n}^{-1/2}X_i^* \e_{i, g_n}^*\|^{2+\dt} \Big| X_0 \right] \nonumber
		\\ \leq & 2^{1+\dt} \left( n^{-\dt/2} h_n^{\dt/2} \sum_{j=1}^{h_n} \ld_j^{-(2+\dt)/2} \right) \left( n^{-1} \sum_{i=1}^n \|X_i\|^{4+2\dt} \right) \|\hat{\beta}_{g_n} - \beta\|^{2+\dt} \label{PBeqBTSvar1_1}
		\\& +  2^{1+\dt} n^{-\dt/2} h_n^{\dt/2} n^{-1} \sum_{i=1}^n \sum_{j=1}^{h_n} \ld_j^{-(2+\dt)/2} |\langle X_i \e_i, \psi_j \rangle|^{2+\dt}. \nonumber
	\end{align}
 based on  
	\begin{align*}
	 	 & \eo^*[\|\Ld_{h_n}^{-1/2}X_i^* \e_{i, g_n}^*\|^{2+\dt}]
		= \left( \sum_{j=1}^{h_n} \ld_j^{-1} \langle X_i\hat{\e}_{i,g_n}, \psi_j \rangle^2 \right)^{(2+\dt) / 2} \\ 
		 & \leq   
		 2^{1+\dt} h_n^{\dt/2} \sum_{j=1}^{h_n} \ld_j^{-(2+\dt)/2} \left(|\langle X_i (\hat{\e}_{i,g_n} - \e_i), \psi_j \rangle|^{2+\dt}  + |\langle X_i \e_i, \psi_j \rangle|^{2+\dt}\right), 
	\end{align*}
	by Jensen's inequality along with
	\begin{align*}
		    \sum_{j=1}^{h_n} \ld_j^{-(2+\dt)/2} |\langle X_i (\hat{\e}_{i,g_n} - \e_i), \psi_j \rangle|^{2+\dt}
		  &=   \sum_{j=1}^{h_n} \ld_j^{-(2+\dt)/2} |\langle \hat{\beta}_{g_n} - \beta,  X_i^{\otimes 2} \psi_j \rangle|^{2+\dt} \\
		 &  \leq      \sum_{j=1}^{h_n} \ld_j^{-(2+\dt)/2}  \|\hat{\beta}_{g_n} - \beta\|^{2+\dt}  \|X_i\|^{4+2\dt} 
	\end{align*}  using
      $\hat{\e}_{i,g_n} - \e_i = - \langle \hat{\beta}_{g_n} - \beta, X_i \rangle$.
	The first term in  \eqref{PBeqBTSvar1_1} converges to zero in probability since $n^{-\dt/2} h_n^{\dt/2} \sum_{j=1}^{h_n} \ld_j^{-(2+\dt)/2} = O(1)$, $\eo[\|X\|^{4+2\dt}] <\infty$, and $\|\hat{\beta}_{g_n} - \beta\| \xrightarrow{\pr} 0$.
	The second term in \eqref{PBeqBTSvar1_1} is bounded as
	\begin{align*}
		\eo \left[ n^{-\dt/2} h_n^{\dt/2} n^{-1} \sum_{i=1}^n \sum_{j=1}^{h_n} \ld_j^{-(2+\dt)/2} |\langle X_i \e_i, \psi_j \rangle|^{2+\dt} \right]
		\leq C n^{-\dt/2} h_n^{(2+\dt)/2}. 
	\end{align*}
	by Condition $\sup_{j \in \N} \ld_j^{-(2+\dt)/2} \eo[|\langle X \e, \psi_j \rangle|^{2+\dt}] < \infty$.
	 Because $n^{-\dt/2}h_n^{(2+\dt)/2} \to 0$ as $n\to\infty$, the second term converges to zero in probability,  verifying \eqref{PBeqBTSvarLyapounov}.
	
	Finally, by combining Slutsky's theorem, Polya's theorem \cite[Theorem 9.1.4]{AL06}, and a subsequence argument (e.g.,~\cite[Theorem~20.5]{bill95}) along with \eqref{PBlemBTSvar2_eq} and \eqref{PBeqBTSvarLyapounov}, we conclude the desired result.
\end{proof}

\begin{proof}[Proof of \autoref{prop_pb_naive_fail}]
	By Propositions~S11-S15~and~Lemma~S56 in the supplement \cite{supp}, 
	the naive bootstrap construction $T_{n,naive}^*(X_0)$ can be written as $T_{n,naive}^*(X_0) = T_{n,\hat{s}}^*(X_0) +A_n^*+ B_n + C_n$, where $T_{n,\hat{s}}^*(X_0)$ is the  PB quantity from \eqref{eq_SnHatStar};
	$A_n^* \equiv A_n^*(X_0)$ is a bootstrap error term that converges to zero in bootstrap probability if $n^{-1/2}h_n^4 (\log h_n)^{7/2} \to 0$; 
	$B_n\equiv B_n(X_0)$ represents a bias-type term that does not depend on the bootstrap sample and  satisfies $\sup_{y \in \R} |\pr(B_n \leq y|X_0) - \Phi(y/\s(\tau))| \xrightarrow{\pr} 0$ with limit variance $\s^2(\tau)$ from \eqref{eq_consbias_limvar}; 
	and $C_n \equiv C_n(X_0)$ is a negligible term that converges to zero if $n^{-1/2} h_n^{9/2} (\log h_n)^6 \to 0$.  By writing $D_n \equiv B_n + C_n$ and
	applying the triangle inequality, we find
	\begin{align*}
		& \left| \sup_{y \in \R} |\pr^*(T_{n,naive}^*(X_0) \leq y | X_0 ) - \pr(T_{n}(X_0) \leq y | X_0 )| - \sup_{y \in \R} |\Phi(y-D_n) - \Phi(y)| \right|
	\end{align*}
	converges to zero in probability, using that $\sup_{y \in \R} |\pr^*(T_{n,\hat{s}}^*(X_0) + A_n^* \leq y  | X_0 ) -  \Phi(y )|\xrightarrow{\pr} 0$, by \autoref{thm_pb}  with Proposition~S11 in \cite{supp},  
	and that $ \sup_{y \in \R} |\pr(T_{n}(X_0) \leq y | X_0 ) - \Phi(y)| \xrightarrow{\pr} 0$ by \autoref{thm_clt}.
	By the continuous mapping theorem/embedding theorem, we then have 
	\begin{align*}
		\Phi (y-D_n ) - \Phi(y) 
		\xrightarrow{\mathsf{d}} \Phi\Big(y+\s(\tau)Z\Big) - \Phi(y),\quad y\in\mathbb{R},
	\end{align*}
  based on  $D_n\equiv B_n + C_n \xrightarrow{\mathsf{d}} -\s(\tau)Z$ for a standard normal  variable $Z$. The convergence in   \autoref{prop_pb_naive_fail} then follows (cf.~\cite{bill95}).
	\end{proof}

	{
	
	\section{Proof of Example~1} \label{app3}
	
	We verify that  \autoref{eg1} conditions ensure   the assumptions for Theorems~\ref{thm_clt}-\ref{thm_pb}. With the given structure, the covariance operator $\Ld \equiv \eo[(X \e)^{\otimes 2}]$ is written as
	\begin{align*}
		\Ld
		&  = \eo[\xi^4] \sum_{j \in \N} \g_j \left( \rho_j^2 \eo[W_1^4]
		+ \sum_{l \neq j} \g_l \rho_l^2 \right) \phi_j^{\otimes 2},
	\end{align*}
	indicating the eigenfunctions $\{\psi_j\}$ of $\Ld$ are given as
	\begin{align*}
		\ld_j
		& = \g_j \eo[\xi^4]\left( \rho_j^2 \eo[W_1^4]
		+ \sum_{l \neq j} \g_l \rho_l^2 \right)
		= \g_j \eo[\xi^4](\ka + 2\rho_j^2) 
	\end{align*}
	 with eigenvalues $\{\ld_j\}$ as $\ld_j  = \phi_j$;  
 note that $\ld_j \asymp \g_j$ as $j\to\infty$. 
	 Condition~(\text{A}1) follows because $\HH = \overline{\mathrm{span}}(\{\phi_j\}_{j=1}^\infty)$ (i.e., the closure of the space spanned by eigenfunctions) and $\ga = \sum_{j=1}^\infty \g_j \phi_j^{\otimes 2}$. Regarding moments conditions for $X$, Condition~(A2) holds because 
			$\eo[(\g_j^{-1/2} \langle X, \phi_j \rangle)^4] = \eo[\xi^4] \eo[W_j^4] < \infty$, while  $\eo[\|X\|^8]<\infty$ holds as  
			$\eo[(\g_j^{-1/2} \langle X, \phi_j \rangle)^8] = \eo[\xi^8] \eo[W_j^8] < \infty$.      
		Conditions (A3)-(A5), all the conditions related the growth rate of $h_n$, and Condition~$B(u)$ with $u>7$ hold by a similar argument to the proof of \cite[Corollary~2]{YDN23RB}.  Conditions~(A6) follows by 
		$h_n^{-1} s_{h_n}(X_0) \to \ka \eo[\xi^4] \xi^2$ here by  
		Lemma~S53 in the supplement.  To verify Condition~(A7),  recall that the 4-th central moment of chi-square distribution with degree of freedom $k$ is $12k^2(k+1)$.
		Then, $\eo[\e^4|X] = 12\|X\|^4 (\|X\|^2+1)$, and hence,
		\begin{align*}
			\eo[\xi_j^4\e^4]
			= \eo[\xi_j^4\eo[\e^4|X]]
			= 12\eo[\xi_j^4\|X\|^4 (\|X\|^2+1)].
		\end{align*}
		As $\|X\|^6 = \left( \sum_{l=1}^\infty \g_l \xi_l^2 \right)^3
		\leq 3 \sum_{l=1}^\infty \g_l^3 \xi_l^6$,
		if we additionally assume $\eo[\xi^{10}]<\infty$, we have
		\begin{align*}
			\sup_{j \in \N} \eo[\xi_j^4\e^4] 
			&  
			\leq \sup_{j \in \N} \sum_{l=1}^\infty \g_l^3 \eo [\xi_j^4 \xi_l^6]
			  \leq C \eo[\xi^{10}] \sum_{l=1}^\infty \g_l^3 \eo[W_1^4] \eo[W_1^6] < \infty;
		\end{align*}
		thus, Condition~(A7) holds as 
		\begin{align*}
			\sup_{j \in \N} \ld_j^{-2} \eo[\langle X\e, \psi_j \rangle^4]
			& = \sup_{j \in \N} \ld_j^{-2} \g_j^2 \eo[\xi_j^4\e^4]
			 = \eo[\xi^4]^{-2}(\ka + \eo[W_1^4]-1)^{-2} \sup_{j \in \N} \eo[\xi_j^4\e^4]  
		\end{align*}
		is finite.  Condition~(A8)  follows by
		$
			\sup_{j \in \N} \g_j^{-1} \langle \Ld \phi_j, \phi_j \rangle
			= \sup_{j \in \N} \eo[\xi^4]\{\ka + \rho_j^2 (\eo[W_1^4]-1)\}
			<\infty$.

}

\end{appendix}

{
 \begin{acks}[Acknowledgments]
 	The authors are grateful to  two anonymous reviewers and an associate editor for their time and  constructive comments that improved the manuscript.
 	Research was partially supported by NSF DMS-2015390.
 \end{acks}
}

%


\bibliographystyle{imsart-number} 
\bibliography{pre_BTS_FLRM_paired_bibfile}       


\end{document}